\documentclass[11pt,a4paper]{article}
\usepackage{fullpage,mathrsfs,graphicx,framed,color,amssymb,amsmath,amsthm}
\usepackage[boxed, noline, ruled, linesnumbered]{algorithm2e}
\usepackage{epsfig, graphicx}
\usepackage{latexsym,amsfonts,amsbsy,amssymb}
\usepackage{amsmath,amsthm}
\usepackage{enumerate}
\usepackage{changes} % for revision
\makeatletter %'@' is now a normal "letter" for TeX

\@addtoreset{equation}{section}
\makeatother %'@' is restored as a "non-letter" character for TeX
\allowdisplaybreaks % For long formula
\newtheorem{theorem}{Theorem}[section]

\newtheorem{remark}{Remark}[section]
\newtheorem{definition}{Definition}[section]
\newtheorem{proposition}{Proposition}[section]
\newtheorem{assumption}{Assumption}

\allowdisplaybreaks % For long formula
\DeclareMathOperator{\essinf}{ess\inf}
\DeclareMathOperator{\esssup}{ess\sup}
\numberwithin{equation}{section}
\begin{document}
\title{Solving High-dimensional Parametric Elliptic Equation Using Tensor Neural
Network\footnote{This work was supported in part
by the National Key Research and Development Program of China (2019YFA0709601 
and 2022YFA1004500), Beijing Natural Science Foundation (Z200003),
National Natural Science Foundations of China (NSFC 11771434 and 12371372), 
Natural Science Foundation of Fujian Province of China (No.2021J01034), 
Fundamental Research Funds for the Central Universities (No.20720220038),
the National Center for Mathematics and Interdisciplinary Science, CAS.}}
\author{Hongtao Chen\footnote{School of Mathematical Sciences 
and Fujian Provincial Key Laboratory of Mathematical Modeling and 
High-Performance Scientific Computing, Xiamen University, Xiamen 361005, China
(chenht@xmu.edu.cn)},\ \ \ Rui Fu\footnote{School of Mathematical Sciences 
and Fujian Provincial Key Laboratory of Mathematical Modeling 
and High-Performance Scientific Computing, Xiamen University, Xiamen 361005, China
(19020211153390@stu.xmu.edu.cn)},\ \ \
Yifan Wang\footnote{LSEC,
Academy of Mathematics and Systems Science,
Chinese Academy of
Sciences, No.55, Zhongguancun Donglu, Beijing 100190, China, and School of
Mathematical Sciences, University of Chinese Academy
of Sciences, Beijing, 100049 (wangyifan@lsec.cc.ac.cn)},\ \ \
and \ \ Hehu Xie\footnote{LSEC,
Academy of Mathematics and Systems Science,
Chinese Academy of
Sciences, No.55, Zhongguancun Donglu, Beijing 100190, China, and School of
Mathematical Sciences, University of Chinese Academy
of Sciences, Beijing, 100049 (hhxie@lsec.cc.ac.cn)} }
\date{}
\maketitle
\begin{abstract}
In this paper, we introduce a tensor neural network based
machine learning method for solving the elliptic partial differential equations
with random coefficients in a bounded physical domain.
With the help of tensor product structure, we can transform the high-dimensional
integrations of tensor neural network functions to one-dimensional integrations
which can be computed with the classical quadrature schemes with high accuracy. 
The complexity of its calculation can be reduced from the exponential scale to 
a polynomial scale. The corresponding machine learning method is designed for
solving high-dimensional parametric elliptic equations.
Some numerical examples are provided to validate the accuracy and
efficiency of the proposed algorithms.

\vskip0.3cm {\bf Keywords.}  Tensor neural network, machine learning,
high-dimensional integration with high accuracy,
parametric elliptic equation with random coefficient.

\vskip0.2cm {\bf AMS subject classifications.} 65N30, 65N25, 65L15, 65B99.
\end{abstract}

\section{Introduction}

Numerical methods for stochastic elliptic partial differential equations (SPDEs)
are vital tools for solving equations that involve both randomness and elliptic operators.
Stochastic elliptic PDEs arise in various fields, including physics, finance, and engineering.
These equations pose challenges due to the interplay between the elliptic operator
and the stochastic component.
Efficient numerical methods for partial differential equations (PDEs)
with random field inputs is a key ingredient
in uncertainty quantification (UQ). %in engineering.
Random field inputs, such as coefficients, loadings and domain,
are spatially inhomogeneous, always lead to the well-known SPDEs \cite{Adler,Vanmarcke}.

There are several numerical methods and techniques that have been developed to tackle SPDEs.
Stochastic finite element method (SFEM) \cite{PapadopoulosGiovanis} 
introduces a stochastic basis expansion to represent
the random coefficients or boundary conditions and solves the resulting deterministic
system of equations using standard finite element method techniques.
SFEM has been extensively studied and applied to a wide range of problems.
Stochastic spectral method employs orthogonal polynomials (e.g., Legendre, Hermite or Laguerre)
to represent the stochastic terms. By applying the Galerkin projection,
the original SPDEs are transformed into a system of deterministic equations,
which can be solved using spectral techniques.
The convergence rate of this method is always faster than that of other methods. 
The generalized polynomial chaos (gPC) method \cite{Xiu} is an extension of extension of the 
spectral method is also widely used for solving SPDEs.  

Monte-Carlo  and quasi-Monte-Carlo sampling strategies \cite{Seydel} of the random inputs,
as well as stochastic collocation methods \cite{Xiu2017}, which have recently 
gained a lot of attention. Monte-Carlo methods simulate the underlying stochastic 
processes and use statistical sampling techniques to obtain approximate solutions.
The idea is to generate a large number of random samples and compute averages
or moments of the solution. For example, Monte-Carlo finite element method,
involves simulating the underlying stochastic processes and solving the deterministic
equations repeatedly for different realizations. By averaging the solutions,
statistical estimates (e.g., mean, variance) of the solution can be obtained.

Due to its universal approximation property, the fully-connected neural network (FNN)
is the most widely used architecture to build the functions for solving high-dimensional
PDEs. There are several types of FNN-based methods such as well-known the deep
Ritz \cite{EYu}, deep Galerkin method \cite{DGM}, PINN \cite{RaissiPerdikarisKarniadakis},
and weak adversarial networks \cite{WAN}
for solving high-dimensional PDEs by designing different loss functions.
Among these methods, the loss functions always include computing high-dimensional
integrations for the functions defined by FNN.
For example,  the loss functions of the deep Ritz method require computing the integrations 
on the high-dimensional domain for the functions which is constructed by FNN.
Direct numerical integration for the high-dimensional functions also
meets the ``curse of dimensionality''. Always, the Monte-Carlo method is
adopted to do the high-dimensional integration with some types of sampling methods
\cite{EYu,HanZhangE}.  Due to the low convergence rate of the Monte-Carlo method,
the solutions obtained by the FNN-based numerical methods are difficult to
obtain high accuracy and stable convergence process. This means that the 
Monte-Carlo method decreases computational work in each forward propagation 
by decreasing the simulation accuracy, efficiency and stability of 
the FNN-based numerical methods for solving high-dimensional PDEs.

Recently, we propose a type of tensor neural network (TNN) and the corresponding
machine learning method
which is designed to solve the high-dimensional problems with high accuracy.
The most important property  of TNN is that the corresponding high-dimensional functions
can be easily integrated with high accuracy and high efficiency. Then the deduced machine
learning method can arrive high accuracy for solving high-dimensional problems.
The reason is that the integration of TNN functions can be separated into one-dimensional
integrations which can be computed by the classical quadrature schemes with high accuracy.
The TNN based machine learning method has already been used to solve high-dimensional
eigenvalue problems and boundary value problems based on the Ritz type of loss
functions \cite{WangJinXie}.
Furthermore,  in \cite{WangXie},  the multi-eigenpairs can also been computed with
machine learning method which is designed by combining the TNN and Rayleigh-Ritz process.
The high accuracy of this type of machine learning methods depends on the fact that the
high-dimensional integrations of functional of TNNs can be computed with high accuracy.

Karhunen-Lo\'{e}ve (KL) expansion approximates the stochastic coefficients
or boundary conditions using a limited number of dominant modes.
This way reduces the dimensionality of the problem
and allows for efficient computations for the deduced deterministic equations which
can then be solved using traditional methods like FEM or FDM.
Unfortunately, the deduced deterministic problems are high-dimensional equations.
Solving these high-dimensional equations by the classical numerical methods always
meets the ``curse of dimensionality''.
In the present paper, we introduce the TNN based machine learning method
for solving the deterministic parametric elliptic equations in a high-(possibly
infinite-) dimensional parameter space, arising as a projection of the SPDE
via a truncated $M$-term KL expansion.
We consider a class of model elliptic problems characterized by
the additive dependence of the equation coefficients
on the multivariate parameter, corresponding to a random field that
is linear in the random variable.

An outline of the paper goes as follows. In Section \ref{Section_SPDE}, the model elliptic
problem with random coefficient and its properties will be presented.
Section \ref{Section_TNN} is devoted to introduce the TNN structure,
its numerical integration method and approximation property.
In Section \ref{Section_TNN_ML}, the TNN based machine learning method
will be proposed to solve the deduced high-dimensional deterministic equation.
Section \ref{Section_Numerical_Examples} gives some numerical examples to validate
the accuracy and efficiency of the proposed TNN based machine learning method.
Some concluding remarks are given in the last section.

\section{Model problem}\label{Section_SPDE}
In this section, we introduce the concerned model problem, the way to
separate stochastic and deterministic variables in the diffusion
coefficient $a(\omega,x)$ with KL expansion, the properties of the
deduced deterministic problem.
The content here is the preparation for designing the TNN based machine learning method for
the high-dimensional parametric problems.

%\subsection{Elliptic problem with stochastic diffusion coefficient}
This paper is mainly concerned with the elliptic problem with stochastic diffusion coefficient
$a(\omega,x)$. We assume the physical domain $D\subset\mathbb{R}^d$, $d\in \mathbb N^+$
is a bounded open set with Lipschitz boundary $\partial D$, and the probability space denoted by $(\Omega,\mathcal{B},P)$,  where $\Omega$ denotes the outcomes,
$\mathcal{B}\subset 2^\Omega$ the sigma algebra of possible events
and $P:\mathcal{B}\rightarrow [0,1]$ a probability measure.
Here, we consider the following problem: Find $u$ such that
\begin{equation}\label{ref:strong}
\Bigg\{
\begin{array}{rcl}
-{\rm div}(a(\omega,x)\nabla u(\omega,x))&=&f(\omega,x), \ \ \ {\rm in}\  D,\\
u(\omega,x)&=&0, \ \ \ \ \ \ \  {\rm on}\ \partial D,
\end{array}
\end{equation}
for $P$-a.e.$\omega\in\Omega$, where $f\in L^2(\Omega\times D)$,
the diffusion coefficient $a\in L^\infty(\Omega\times D)$
and the solution $u(\omega,x)$ are random fields in the spatial domain $D$.

In order to guarantee the existence, uniqueness and well-posedness of the problem,
following \cite{AbdulleBarthSchwab},
we make the following assumption on the diffusion coefficient $a(\omega,x)$.
\begin{assumption}\label{Assumption_1}
Let $a\in L^\infty(\Omega\times D)$ be strictly positive, with lower and upper bound
\begin{eqnarray}
P\left\{w\in\Omega:a_{\rm min}\leq \operatornamewithlimits{\essinf}_{x\in D} a(\omega,x) \land \operatornamewithlimits{\esssup}_{x\in D}a(\omega,x)\leq a_{\rm max}\right\}=1,
\end{eqnarray}
where the essential infimum and supremum are taken with respect to the Lebesgue
measure in $D\subset\mathbb{R}^d$.
\end{assumption}

%\subsection{Variational Formulation}
In order to formulate the variational formulation of the problem, 
we introduce the tensor Hilbert space
\begin{eqnarray}
L^2_P(\Omega;H^1_0(D)):=L^2_P(\Omega)\otimes H^1_0(D),
\end{eqnarray}
 where
\begin{eqnarray*}
L^2_P(\Omega):=\left\{\xi(\omega)\Big\vert\int_\Omega\xi^2(\omega)dP(\omega)<\infty\right\}.
\end{eqnarray*}

Then we multiply (\ref{ref:strong}) with a test function $v\in L^2_P(\Omega;H^1_0(D))$,
integrate by parts over the physical domain $D$ and take expectations on both sides,
yielding the variational formulation for (\ref{ref:strong}):
Find $u\in L^2_P(\Omega;H^1_0(D))$ such that for all $v\in L^2_P(\Omega;H^1_0(D))$
\begin{eqnarray}\label{ref:weak}
b(u,v):=\mathbb{E}\bigg[\displaystyle\int_D a(\omega,x)\nabla u(\omega,x)\cdot\nabla v(\omega,x)dx\bigg]=\mathbb{E}\bigg[\displaystyle\int_D f(\omega,x)v(\omega,x)dx\bigg].
 \end{eqnarray}
Under Assumption \ref{Assumption_1}, the bilinear form $b(\cdot,\cdot)$ is
coercive and continuous. Hence, by using Lax-Milgram lemma, the equation (\ref{ref:weak})
is uniquely solvable for every $f\in L^2_P(\Omega;H^{-1}(D))$.

%\subsection{Separation of stochastic and deterministic variables}
Here, the KL  expansion is used to separate stochastic
and deterministic variables in the diffusion coefficient $a(\omega,x)$.
First we assume the information about the diffusion coefficient
is known including its mean field
\begin{eqnarray}\label{Mean_Definition}
\mathbb{E}_a(x)= \displaystyle\int_\Omega a(\omega,x)dP(\omega),\enspace x\in D,
\end{eqnarray}
and covariance
\begin{eqnarray}\label{covariance}
C_a(x,x')=\displaystyle\int_\Omega \left(a(\omega,x)
-\mathbb{E}_a(x)\right)\left(a(\omega,x')
-\mathbb{E}_a(x')\right)dP(\omega),\enspace x,x'\in D.
\end{eqnarray}
In this paper, we assume  $C_a(x,x')$ is an admissible covariance function
which is defined as follows.
\begin{definition}\label{SPD_Definition}
A covariance function $C_a(x,x')\in L^2(D\times D)$ given by (\ref{covariance})
is said to be admissible if it is symmetric and positive definite in the following sense
\begin{eqnarray}
0\leq\sum\limits_{k=1}^n\sum\limits_{j=1}^n a_k C_a(x_k,x_j)\overline{a}_j,
\enspace\enspace\forall x_k,x_j\in D,a_k,a_j\in\mathbb{C}.
\end{eqnarray}
\end{definition}
Based on the random permeability $a(\omega,x)\in L^2(\Omega\times D)$, we can define
the covariance operator as
\begin{eqnarray}
\mathcal{C}_a:L^2(D)\rightarrow L^2(D),\enspace(\mathcal{C}_a u)(x)
:=\displaystyle\int_D C_a(x,x')u(x')dx'.
\end{eqnarray}
Given an admissible covariance function $C_a(x,x')$ in the sense of Definition \ref{SPD_Definition},
the associated covariance operator $\mathcal{C}_a$ is a symmetric,
non-negative and compact integral operator from $L^2(D)$ to $L^2(D)$.
Therefore, it has a countable sequence of eigenpairs $(\lambda_m,\phi_m)_{m\geq1}$
\begin{eqnarray*}
\mathcal{C}_a\varphi_m=\lambda_m\varphi_m,\enspace\enspace m=1,2\ldots,
\end{eqnarray*}
where the sequence (real and positive) of KL-eigenvalues $\lambda_m$
is enumerated with decreasing magnitude and is either finite or tends to zero
as $m\rightarrow\infty$.
The corresponding KL eigenfunctions $\varphi_m(x)$ are assumed to be $L^2(D)$-orthonormal, i.e.
\begin{eqnarray*}
\displaystyle\int_D\varphi_m(x)\varphi_n(x)d x=\delta_{mn}, \quad m,n=1,2, \ldots.
\end{eqnarray*}
Then we have the following expansion for the diffusion coefficient $a(\omega, x)$.
\begin{definition}
The KL expansion of a random field $a(\omega,x)$ with finite
mean (\ref{Mean_Definition}) and  admissible covariance (\ref{covariance}) is given by
\begin{eqnarray}\label{ref:KL_expansion}
a(\omega,x)=\mathbb{E}_a(x)+\sum\limits_{m\geq 1}\sqrt{\lambda_m}\varphi_m(x)Y_m(\omega).
\end{eqnarray}
The family of random variables $(Y_m)_{m\geq1}$ is given by
\begin{eqnarray}\label{ref:variables}
Y_m(\omega)=\displaystyle\frac{1}{\sqrt{\lambda_m}}\displaystyle\int_D \left(a(\omega,x)-\mathbb{E}_a(x)\right)\varphi_m(x)dx:\Omega\rightarrow\mathbb{R}.
\end{eqnarray}
\end{definition}
It should be pointed out that $Y_m(\omega)$ are centered with unit variance 
and pairwise uncorrelated.

%\subsection{Parametric deterministic problem}
In order to parametrize the stochastic input, we make the following assumption
on the random variables $Y_m$ in the KL-representation $(\ref{ref:KL_expansion})$ of $a(\omega,x)$.
\begin{assumption}\label{Assumption_2}
The following assumptions hold
\begin{enumerate}
\item [(1)]
The family $(Y_m)_{m\geq1}:\Omega\rightarrow\mathbb{R}$ is independent.

\item [(2)] The KL-expansion $(\ref{ref:KL_expansion})$ of the input data is $finite$,
i.e. there exists $\bar{M}<\infty$ such that $Y_m=0$ for all $m>\bar{M}$.

\item [(3)] Each $Y_m(\omega)$ in (\ref{ref:KL_expansion}), (\ref{ref:variables})
is associated a complete probability space $(\Omega_m, \Sigma_m, P_m)$,
$m\in \mathbb{N}$ with the following properties:
\begin{itemize}
\item [(a)] The range $\Gamma_m:={\rm Ran}(Y_m)\subseteq\mathbb{R}$ of $Y_m$
is compact and is scaled to $[-1,1]$ for all $m$.

\item[(b)] The probability measure $P_m$ deduces a probability density function
$\rho_m:\Gamma_m\rightarrow[0,\infty]$ such that $dP_m(\omega)=\rho_m(y_m)dy_m$,
$m\in\mathbb{N}$, $y_m\in\Gamma_m$.

\item[(c)] The sigma algebras $\Sigma_m$ are subsets of the Borel sets on
the interval $\Gamma_m$, i.e. $\Sigma_m\subseteq \mathcal{B}(\Gamma_m)$.
\end{itemize}
\end{enumerate}
\end{assumption}

In order to numerically treat the KL-expansion (\ref{ref:KL_expansion}),
we define the $M$-term truncated KL-expansion as follows
\begin{eqnarray}\label{ref:truncated_KL_expansion}
a_M(\omega,x)=\mathbb{E}_a(x)+\sum\limits_{m=1}^M\sqrt{\lambda_m}\varphi_m(x)Y_m(\omega).
\end{eqnarray}
Replacing the random input $a(\omega,x)$ in (\ref{ref:strong})
by its $M$-term truncated KL-expansion $a_M(\omega,x)$
in (\ref{ref:truncated_KL_expansion}) leads to a $M$-dimensional parametric,
deterministic variational formulation. 
In addition, we further assume that the right-hand side $f(\omega,x)$ can be represented as 
a finite-term K-L expansion as follows
\begin{eqnarray}\label{ref:truncated_KL_expansion}
f(\omega,x)=\mathbb{E}_f(x)+\sum\limits_{m=1}^M\sqrt{\mu_m}\psi_m(x)Y_m(\omega).
\end{eqnarray}
\begin{remark}
This setting is referred to \cite{BabuskaNobileTempone}.
Since it is usual to have $f$ and $a$ independent, because the loads and the 
material properties are seldom related.
In such a situation we have $a(Y_{1}(\omega),\cdots,Y_M(\omega),x)=a(Y_{a,1}(\omega),\cdots,Y_{a,M_a}(\omega),x)$
and $f(Y_{1}(\omega),\cdots,Y_M(\omega),x)=f(Y_{f,1}(\omega),\cdots,Y_{f,M_f}(\omega),x)$,
with random vectors $(Y_{a,1},\cdots,Y_{a,M_a})$ and $(Y_{f,1},\cdots,Y_{f,M_f})$ independent, $M=M_a+M_f$.
\end{remark}

In order to define the variational form for the stochastic problem (\ref{ref:strong}),
we introduce the following function spaces: For
\begin{eqnarray*}
\Gamma:=\Gamma_1\times\Gamma_2\times\cdots,\enspace \ \ \
y=(y_1,y_2,\cdots)\in\Gamma,
\end{eqnarray*}
we denote the space of all mappings $v:\Gamma\rightarrow\mathbb{R}$
which are square integrable with respect to the measure $\rho(dy)=\rho(y)dy$
as $L^2_{\rho}(\Gamma)$, where $\rho(dy)$ is a probability measure on $\Gamma$.
If $H$ denotes a separable Hilbert space with norm $\Vert\cdot\Vert_H$,
we denote the Bochner space of functions $v:\Gamma\rightarrow H$
such that $\Vert v(y,\cdot)\Vert_H:y\rightarrow
\mathbb{R}$ belongs to $L^2_\rho(\Gamma)$. In this paper, the solution is
required in the space
\begin{eqnarray*}\label{ref:space1}
L^2_\rho(H_0^1):=L^2_\rho(\Gamma;H_0^1(D))\simeq L^2_\rho(\Gamma)\otimes H_0^1(D),
\end{eqnarray*}
where $\otimes$ denotes the tensor product between separable Hilbert spaces.

By Assumption \ref{Assumption_2}(1), the probability density $\rho(y)$ is separable, i.e.
\begin{eqnarray}\label{Density}
\rho(y)=\prod\limits_{m\geq 1}\rho_m(y_m)
\end{eqnarray}
with $\rho_m(y_m)$ satisfying Assumption \ref{Assumption_2}(3)$(b)$.
Based on the above assumptions, we can obtain the parametric,
deterministic formulation of (\ref{ref:strong}) in variational form:
Find $u_M\in L^2_\rho(H^1_0)$ such that
\begin{eqnarray}\label{Weak_Problem}
b_M(u_M,v)=\ell(v),\enspace\enspace\forall v\in L^2_\rho(H^1_0),
\end{eqnarray}
with
\begin{eqnarray}
b_M(u_M,v)&=& \mathbb{E} \Bigg[\displaystyle
\int_Da_M(y,x)\nabla u_M(y,x)\cdot\nabla v(y,x)dx\Bigg]\nonumber\\
&=&\displaystyle\int_\Gamma\displaystyle 
\int_Da_M(y,x)\nabla u_M(y,x)\cdot\nabla v(y,x)\rho(y)dxdy,
\end{eqnarray}
and
\begin{eqnarray}
\ell(v)=\mathbb{E} \Bigg[\displaystyle\int_D f(y,x)v(y,x)dx\Bigg]
=\displaystyle\int_\Gamma\displaystyle \int_D f(y,x)v(y,x)\rho(y)dxdy.
\end{eqnarray}
Here, based on (\ref{ref:truncated_KL_expansion}), the truncated KL-expansion has following
expansions.
\begin{eqnarray}
a_M(y,x)=a_0(x)+\sum \limits_{m=1}^M y_m\psi_m(x),\enspace a_0=\mathbb{E}_a,\enspace \psi_m(x):=\lambda_m^{1/2}\varphi_m(x).
\end{eqnarray}

%\subsection{Regularity}
Now, we come to introduce some regularity results of the concerned problem (\ref{Weak_Problem}).
In \cite{TodorSchwab}, it is proved that the weak solution $u(y, \cdot) \in H_0^1(D)$
is analytic as a function of
$$
y \mapsto u(y, \cdot) \in H_0^1(D) \quad \text { from} \  \Gamma \text { to } H_0^1(D).
$$
The precise quantitative analysis on the size of the domain of
analyticity is based on the decay rate of the coefficients
in (\ref{ref:truncated_KL_expansion}).
Here we consider two basic cases of coefficient decay:
\begin{enumerate}
\item  Exponential decay (see \cite{TodorSchwab}):
\begin{eqnarray}\label{Exponential_decay}
\rho_m:=\left\|a_m\right\|_{L^{\infty}(D)} \leq C_0 \exp \left(-C_1 m^{1 / d}\right),
\quad \forall m \in \mathbb{N}_{+}.
\end{eqnarray}
Note that the sequence $\rho=\left(\rho_m\right)_{m=1}^{\infty}$ in
(\ref{Density}) belongs to $\ell^p(\mathbb{N})$ for every $p>0$.
\item  Algebraic decay (see, e.g., \cite{BieriAndreevSchwab,CohenDeVoreSchwa}):
there is a constant $s>0$ such that
the following bounds hold
\begin{eqnarray}\label{Algebraic_decay}
\rho_m:=\left\|a_m\right\|_{L^{\infty}(D)} \leq C_0 m^{-s / d}, \quad \forall m \in \mathbb{N}_{+}.
\end{eqnarray}
Note that the sequence $\rho=\left(\rho_m\right)_{m=1}^{\infty}$ in (\ref{Density})
belongs to $\ell^p(\mathbb{N})$ for every $p>d_0 / s$.
\end{enumerate}
Based on the coefficient decay condition (\ref{Exponential_decay}), it has been proved
in \cite{TodorSchwab} that the domain of analyticity of the solution $u$ of (\ref{ref:strong})
as a function of $y_m$ increases exponentially in size as $m \nearrow \infty$.
The following proposition gives the explicit bounds on all derivatives of $u$ with respect to
$y_m$ which are suitable for either case 1 or 2.
\begin{proposition} (see \cite{TodorSchwab}).
If the exponential decay condition (\ref{Exponential_decay}) holds,
$u=u_M(y, \cdot)$ is the exact solution of (\ref{ref:strong}) and $M$ is large enough,
then the following estimate
\begin{eqnarray}\label{Regularity_Result}
\left\|\partial_y^\alpha u_M(y, \cdot)\right\|_{H_0^1(D)}
\leq\left(C_a^{|\alpha|} \cdot|\alpha|! \cdot \prod_{m=1}^M \rho_m^{\alpha_m}\right)
\end{eqnarray}
for all $y \in \Gamma$ and $\alpha \in \mathbb{N}_0^M$ where $\mathbb N_0$ denotes the 
set of all non-negative integers.
\end{proposition}

\section{Tensor neural network and its quadrature scheme}\label{Section_TNN}
In this section, we introduce the TNN structure and the quadrature scheme for the high-dimensional
TNN functions. %\magenta{(There include)}
Also included in this section are a discussion of the approximation properties, some techniques to
improve the numerical stability, and the complexity estimate of the 
high-dimensional integrations of the TNN functions.

%-----------------------------------------------------------------------------------------------------
\subsection{Tensor neural network architecture}
This subsection is devoted to introducing the TNN structure and some techniques to 
improve the stability of the corresponding machine learning methods.
The approximation properties of TNN functions have been discussed and 
investigated in \cite{WangJinXie}.
In order to express clearly and facilitate the construction of the TNN method for solving
high-dimensional PDEs, here we will also introduce some important definitions and properties.

\begin{figure}[htb!]
\centering
\includegraphics[width=16cm,height=12cm]{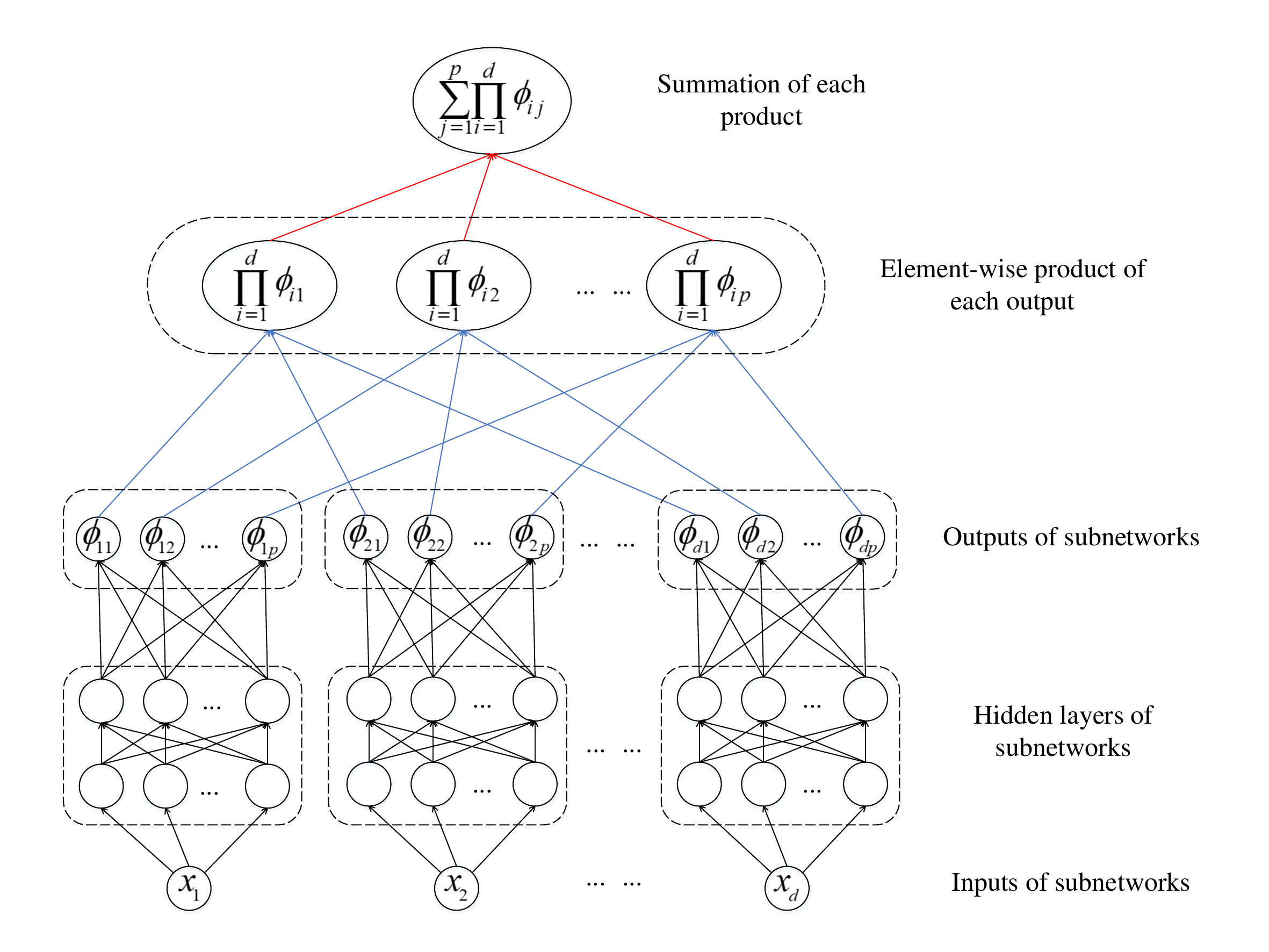}
\caption{Architecture of TNN. Black arrows mean linear transformation
(or affine transformation). Each ending node of blue arrows is obtained by taking the
scalar multiplication of all starting nodes of blue arrows that end in this ending node.
The final output of TNN is derived from the summation of all 
starting nodes of red arrows.}\label{TNNstructure}
\end{figure}

TNN is built by the tensor products of one dimension functions which come from
$d$ subnetworks with one-dimensional input and multi-dimensional output, where
$d$ denotes the spatial dimension of the concerned problems which will be solved by the
machine learning method in this paper.
For each $i=1,2,\cdots,d$, we use $\Phi_i(x_i;\theta_i)
=(\phi_{i,1}(x_i;\theta_i),\phi_{i,2}(x_i;\theta_i),\cdots,\phi_{i,p}(x_i;\theta_i))$
to denote a subnetwork that maps a set $\Omega_i\subset\mathbb R$ to $\mathbb R^p$,
where $\Omega_i,i=1,\cdots,d,$ can be a bounded interval $(a_i,b_i)$, 
the whole line $(-\infty,+\infty)$
or the half line $(a_i,+\infty)$.
The number of layers and neurons in each layer, the selections of activation 
functions and other hyperparameters can be different in different subnetworks. 
TNN consists of $p$ correlated rank-one functions,
which are composed of the multiplication of $d$ one-dimensional input 
functions in different directions.
Figure \ref{TNNstructure} shows the corresponding architecture of TNN.

In order to improve the numerical stability, we normalize each $\phi_{i,j}(x_i)$
and use the following normalized-TNN structure:
\begin{eqnarray}\label{def_TNN_normed}
\Psi(x;c,\theta)=\sum_{j=1}^pc_j\widehat\phi_{1,j}(x_1;\theta_1)\widehat\phi_{2,j}(x_2;\theta_2)
\cdots\widehat\phi_{d,j}(x_d;\theta_d)
=\sum_{j=1}^pc_j\prod_{i=1}^d\widehat\phi_{i,j}(x_i;\theta_i),
\end{eqnarray}
where each $c_j$ is a scaling parameter which describes the length of each rank-one function, $c=\{c_j\}_{j=1}^{p}$ is a set of trainable parameters, $\{c,\theta\}=\{c,\theta_1,\cdots,\theta_d\}$
denotes all parameters of the whole architecture.
For $i=1,\cdots,d,j=1,\cdots,p$, $\widehat\phi_{i,j}(x_i,\theta_i)$ 
is a $L^2$-normalized function as follows:
\begin{eqnarray}\label{eq_phi_normed}
\widehat\phi_{i,j}(x_i,\theta_i)
=\frac{\phi_{i,j}(x_i,\theta_i)}{\|\phi_{i,j}(x_i,\theta_i)\|_{L^2(\Omega_i)}}.
\end{eqnarray}
For simplicity of notation, $\phi_{i,j}(x_i,\theta_i)$ denotes 
the normalized function in the following parts.

The TNN architecture (\ref{def_TNN_normed}) and the architecture defined in 
\cite{WangJinXie} are mathematically equivalent, 
but (\ref{def_TNN_normed}) has better numerical stability during the training process.
From Figure \ref{TNNstructure} and numerical tests, we can find the parameters for 
each rank of TNN are correlated by the FNN, which guarantee the stability 
of the TNN-based machine learning methods. This is also an important difference 
from the tensor finite element methods.

In order to show the reasonableness of TNN, we now introduce the approximation property 
from \cite{WangJinXie}.
Since there exists the isomorphism relation between $H^m(\Omega_1\times\cdots\times\Omega_d)$
and the tensor product space $H^m(\Omega_1)\otimes\cdots\otimes H^m(\Omega_d)$,
the process of approximating the function $f(x)\in H^m(\Omega_1\times\cdots\times\Omega_d)$
by the TNN defined as (\ref{def_TNN_normed}) can be regarded as searching for 
a correlated CP decomposition structure
to approximate $f(x)$ in the space $H^m(\Omega_1)\otimes\cdots\otimes H^m(\Omega_d)$
with the rank being not greater than $p$.
The following approximation result to the functions in the space 
$H^m(\Omega_1\times\cdots\times\Omega_d)$ under the sense of $H^m$-norm 
is proved in \cite{WangJinXie}.
\begin{theorem}\cite[Theorem 1]{WangJinXie}\label{theorem_approximation}
Assume that each $\Omega_i$ is an interval in $\mathbb R$ for $i=1, \cdots, d$, $\Omega=\Omega_1\times\cdots\times\Omega_d$,
and the function $f(x)\in H^m(\Omega)$. Then for any tolerance $\varepsilon>0$, there exist a
positive integer $p$ and the corresponding TNN defined by (\ref{def_TNN_normed})
such that the following approximation property holds
\begin{equation}\label{eq:L2_app}
\|f(x)-\Psi(x;\theta)\|_{H^m(\Omega)}<\varepsilon.
\end{equation}
\end{theorem}

Although there is no general result to give the relationship between the 
hyperparameter $p$ and error bounds,
we also provide an estimate of the rank $p$ under a smoothness assumption. 
For easy understanding, we focus on the periodic setting with 
$I^d=I\times I\times\cdots\times I=[0,2\pi]^d$ and
the approximation property of TNN to the functions in the 
linear space which is defined by Fourier basises.
Note that similar approximation results of TNN can be extended to the non-periodic functions.
For each variable $x_i\in[0,2\pi]$, let us define the one-dimensional
Fourier basis $\{\varphi_{k_i}(x_i):= \frac{1}{\sqrt{2\pi}}e^{{\rm i}k_ix_i},k_i\in\mathbb Z\}$
and classify functions via the decay of their Fourier coefficients.
Further we denote multi-index $k=(k_1,\cdots,k_d)\in\mathbb Z^d$ and $x=(x_1,\cdots,x_d)\in I^d$.
Then the $d$-dimensional Fourier basis can be built in the tensor product way
\begin{eqnarray}
\varphi_k(x):=\prod_{i=1}^d\varphi_{k_i}(x_i)=(2\pi)^{-d/2}e^{{\rm i}k\cdot x}.
\end{eqnarray}
We denote
\begin{eqnarray}
\lambda_{\rm mix}(k):=\prod_{i=1}^d\left(1+|k_i|\right)\ \ \ {\rm and}\ \ \ 
\lambda_{\rm iso}(k):= 1+\sum_{i=1}^d|k_i|.
\end{eqnarray}
Now, for $-\infty<t,\ell<\infty$, we define the space $H_{\rm mix}^{t,\ell}(I^d)$ 
as follows (cf. \cite{GriebelHamaekers,Knapek})
\begin{eqnarray}
H_{\rm mix}^{t,\ell}(I^d)=\left\{u(x)=\sum_{k\in\mathbb Z^d}c_k\varphi_k(x): \|u\|_{H_{\rm mix}^{t,\ell}(I^d)}=\left(\sum_{k\in\mathbb Z^d}\lambda_{\rm mix}(k)^{2t}\cdot\lambda_{\rm iso}(k)^{2\ell}\cdot|c_k|^2\right)^{1/2}<\infty \right\}.
\end{eqnarray}
Note that the parameter $\ell$ governs the isotropic smoothness, 
whereas $t$ governs the mixed smoothness.
The space $H_{\rm mix}^{t,\ell}(I^d)$ gives a quite flexible framework 
for the study of problems in Sobolev spaces.
See \cite{GriebelHamaekers,GriebelKnapek,Knapek} for more information on the space 
$H_{\rm mix}^{t,\ell}(I^d)$.

Thus, \cite{WangJinXie} gives the following comprehensive error estimate for TNN.
\begin{theorem}\label{theorem_aprrox_rate}
Assume function $f(x)\in H_{\rm mix}^{t,\ell}(I^d)$, $t>0$ and $m>\ell$. 
Then there exists a TNN $\Psi(x;\theta)$ defined by (\ref{def_TNN_normed}) 
such that the following approximation property holds
\begin{eqnarray}
\|f(x)-\Psi(x;\theta)\|_{H^m(I^d)}\leq C(d)\cdot p^{-(\ell-m+t)}\cdot\|u\|_{H_{\rm mix}^{t,\ell}(I^d)},
\end{eqnarray}
where $C(d)\leq c\cdot d^2\cdot0.97515^d$ and $c$ is independent of $d$.
And each subnetwork of TNN is a FNN which is built by
using $\sin(x)$ as the activation function and one hidden layer with $2p$ neurons,
see Figure \ref{TNNstructure}.
\end{theorem}

The TNN-based machine learning method in this paper will adaptively 
select $p$ rank-one functions by training process.
From the approximation result in Theorem \ref{theorem_aprrox_rate}, 
when the target function belongs to $H_{\rm mix}^{t,\ell}(\Omega)$,  
there exists a TNN with $p\sim\mathcal O(\varepsilon^{-(m-\ell-t)})$ 
such that the accuracy is $\varepsilon$.
For more details about the approximation properties of TNN, please refer to \cite{WangJinXie}.

\subsection{Quadrature scheme for TNN}\label{section_quad}
In this section, we introduce the quadrature scheme for computing 
the high-dimensional integrations of the
high-dimensional TNN functions.
Due to the low-rank of TNN, the efficient and accurate quadrature scheme can be designed for the
TNN-related high-dimensional integrations which are included in the loss
functions for machine learning methods.
The detailed numerical integration scheme for TNN functions
with the polynomial scale computational complexity of the dimension will be designed here.
When calculating the high-dimensional integrations in the loss functions,
we only need to calculate one-dimensional integrations as well as their product.

The method to compute the numerical integrations of polynomial composite
functions of TNN and their derivatives has already designed in \cite{WangJinXie}.
For the sake of completeness, we also introduce the quadrature scheme here.
For more information, please refer to \cite{WangJinXie}.

In order to describe the differentiation of TNNs and their composition with a polynomial,
we define the following index sets
\begin{eqnarray*}
\mathcal B&:=&\left\{\beta=(\beta_1,\cdots,\beta_d)\in\mathbb N_0^d\ 
\Big|\ |\beta|:=\sum_{i=1}^d\beta_i\leq m \right\},\\
\mathcal A&:=&\left\{\alpha=(\alpha_\beta)_{\beta\in\mathcal B}\in\mathbb N_0^{|\mathcal B|}\
\Big|\ |\alpha|:=\sum_{\beta\in\mathcal B}\alpha_\beta\leq k \right\},
\end{eqnarray*}
where $m$ and $k$ are two positive integers,
$|\mathcal B|$ and $|\mathcal A|$ denote the cardinal 
numbers of $\mathcal B$ and $\mathcal A$, respectively.
It is easy to know that
$$|\mathcal B|=\sum_{j=0}^m\binom{j+d-1}{j}, \ \ \ \ 
|\mathcal A|=\sum_{j=0}^k\binom{j+|\mathcal B|-1}{j},$$
and the scales of magnitudes of $|\mathcal B|$ and $|\mathcal A|$ are
$\mathcal O\big((d+m)^m\big)$
and $\mathcal O\big(((d+m)^m+k)^k\big)$, respectively.

In the following parts of this paper, the parameters in (\ref{def_TNN_normed}) 
will be omitted for brevity without confusion. The activation function of TNN is selected
to be smooth enough such  that $\Psi(x)$ has partial derivatives up to order $m$.
In this section, we assume the integrand, denoted by $F(x)$,
includes the $k$-degree complete polynomial of $d$-dimensional TNN and
its partial derivatives up to order $m$ that can be described as follows
\begin{eqnarray}\label{def_F(x)}
F(x)=\sum_{\alpha\in\mathcal A}A_{\alpha}(x)\prod_{\beta\in\mathcal B}
\left(\frac{\partial^{|\beta|}\Psi(x)}{\partial x_1^{\beta_1}
\cdots\partial x_d^{\beta_d}}\right)^{\alpha_\beta},
\end{eqnarray}
where the coefficient $A_\alpha(x)$ has the decomposition form
with rank less than $q$ in the tensor product space $L^2(\Omega_1)\otimes\cdots\otimes L^2(\Omega_d)$
\begin{eqnarray}\label{def_A_alpha}
A_\alpha(x)=\sum_{\ell=1}^q B_{1,\ell,\alpha}(x_1)B_{2,\ell,\alpha}(x_2)\cdots B_{d,\ell,\alpha}(x_d),
\end{eqnarray}
where $B_{i,\ell,\alpha}(x_i)$ denotes the one-dimensional function in 
$L^2(\Omega_i)$ for $i=1, \cdots, d$
and  $\ell=1,\cdots, q$.
Thanks to the tensor product structure of $F(x)$, the
high-dimensional multiple integral $\int_\Omega F(x)dx$ for the TNN
can be decomposed into the sum of products of one-dimensional integrations
which can be computed by common numerical quadrature.
In order to describe  the  decomposition, for each
$\alpha = (\alpha_1,\cdots,\alpha_{|\mathcal B|})\in\mathcal A$, we  define the following set
\begin{eqnarray*}
\mathcal B_\alpha:=\Big\{\beta=\left(\beta_1,\cdots,\beta_d\right)\in\mathcal B
\ \big|\ \alpha_\beta\geq 1\Big\}.
\end{eqnarray*}
By the definition of the index set $\mathcal A$,  the cardinal of $\mathcal B_\alpha$ 
has the estimate $|\mathcal B_\alpha|\leq k$ for any $\alpha\in\mathcal A$.

With the help of TNN structure  (\ref{def_TNN_normed}) of  $\Psi(x)$, 
the products in (\ref{def_F(x)}) can be further decomposed as
\begin{eqnarray}\label{eq_decomposition_prod}
&&\prod_{\beta\in\mathcal B_\alpha}
\left(\frac{\partial^{|\beta|}\Psi(x)}{\partial x_1^{\beta_1}
\cdots\partial x_d^{\beta_d}}\right)^{\alpha_\beta}\nonumber\\
&&=\sum_{\substack{\beta\in\mathcal B_\alpha, \ell=1,\cdots,\alpha_\beta, \\ 
1\leq j_{\beta,\ell}\leq p}}
\left(\prod_{\beta\in\mathcal B_\alpha}\prod_{\ell=1}^{\alpha_\beta}
\frac{\partial^{\beta_1}\phi_{1,j_{\beta,\ell}}(x_1)}{\partial x_1^{\beta_1}}\right)\cdots
\left(\prod_{\beta\in\mathcal B_\alpha}\prod_{\ell=1}^{\alpha_\beta}
\frac{\partial^{\beta_d}\phi_{d,j_{\beta,\ell}}(x_d)}{\partial x_d^{\beta_d}}\right).\ \ \ \ \
\end{eqnarray}
With the help of expansion (\ref{eq_decomposition_prod}),  
we can give the following expansion for $F(x)$
\begin{eqnarray}\label{Expansion_Fx}
F(x)&=&\sum_{\alpha\in\mathcal A}\sum_{\ell=1}^q\sum_{\substack{\beta\in\mathcal B_\alpha,
\ell=1,\cdots,\alpha_\beta,\\ 1\leq j_{\beta,\ell}\leq p}}
\left(B_{1,\ell,\alpha}(x_1)\prod_{\beta\in\mathcal B_\alpha}\prod_{\ell=1}^{\alpha_\beta}
\frac{\partial^{\beta_1}\phi_{1,j_{\beta,\ell}}(x_1)}{\partial x_1^{\beta_1}}\right)\nonumber\\
&&\ \quad\quad \ \cdots\left(B_{d,\ell,\alpha}(x_d)
\prod_{\beta\in\mathcal B_\alpha}\prod_{\ell=1}^{\alpha_\beta}
\frac{\partial^{\beta_d}\phi_{d,j_{\beta,\ell}}(x_d)}{\partial x_d^{\beta_d}}\right).
\end{eqnarray}
Based on the decomposition (\ref{Expansion_Fx}), 
the high-dimensional integration $\int_\Omega F(x)dx$
can be decomposed into the following one-dimensional integrations
\begin{eqnarray}\label{Integration_Expansions}
\int_\Omega F(x)dx &=&\sum_{\alpha\in\mathcal A}\sum_{\ell=1}^q\sum_{\substack{\beta\in\mathcal B_\alpha, \ell=1,\cdots,\alpha_\beta,\\ 1\leq j_{\beta,\ell}\leq p}}
\int_{\Omega_1}\left(B_{1,\ell,\alpha}(x_1)\prod_{\beta\in\mathcal B_\alpha}\prod_{\ell=1}^{\alpha_\beta}
\frac{\partial^{\beta_1}\phi_{1,j_{\beta,\ell}}(x_1)}{\partial x_1^{\beta_1}}\right)dx_1\nonumber\\
&&\ \quad\quad \ \cdots\int_{\Omega_d}\left(B_{d,\ell,\alpha}(x_d)
\prod_{\beta\in\mathcal B_\alpha}\prod_{\ell=1}^{\alpha_\beta}
\frac{\partial^{\beta_d}\phi_{d,j_{\beta,\ell}}(x_d)}{\partial x_d^{\beta_d}}\right)dx_d.
\end{eqnarray}
The expansion (\ref{Integration_Expansions}) gives the hint to design the efficient numerical
scheme to compute the high-dimensional integration of the TNN function $F(x)$.
Without loss of generality, for the $i$-th one-dimensional domain $\Omega_i$,
we choose $N_i$ Gauss points $\{x_i^{(n_i)}\}_{n_i=1}^{N_i}$
and the corresponding weights $\{w_i^{(n_i)}\}_{n_i=1}^{N_i}$ to compute 
the one-dimensional integrations included in  (\ref{Integration_Expansions}), where $i=1, \cdots, d$.
Then we can give the following splitting numerical integration scheme for $\int_\Omega F(x)dx$:
\begin{eqnarray}\label{eq_I_tensor_form}
\int_\Omega F(x)dx&\approx&
\sum_{\alpha\in\mathcal A}\sum_{\ell=1}^q\sum_{\substack{\beta\in\mathcal B_\alpha,
\ell=1,\cdots,\alpha_\beta,\\ 1\leq j_{\beta,\ell}\leq p}}
\left(\sum_{n_1=1}^{N_1}w_1^{(n_1)}B_{1,\ell,\alpha}(x_1^{(n_1)})\prod_{\beta\in\mathcal B_\alpha}\prod_{\ell=1}^{\alpha_\beta}
\frac{\partial^{\beta_1}\phi_{1,j_{\beta,\ell}}(x_1^{(n_1)})}{\partial x_1^{\beta_1}}\right)\nonumber\\
&&\ \ \quad \ \cdots\left(\sum_{n_d=1}^{N_d}w_d^{(n_d)}B_{d,\ell,\alpha}(x_d^{(n_d)})
\prod_{\beta\in\mathcal B_\alpha}\prod_{\ell=1}^{\alpha_\beta}
\frac{\partial^{\beta_d}\phi_{d,j_{\beta,\ell}}(x_d^{(n_d)})}{\partial x_d^{\beta_d}}\right).
\end{eqnarray}
The simplicity of the one-dimensional integration makes the scheme  (\ref{eq_I_tensor_form})
reduce the computational work of the high-dimensional integration of the $d$-dimensional
function $F(x)$ to the polynomial scale of dimension $d$.
Let us denote $N=\max\{N_1,\cdots,N_d\}$ and $\underline N = \min\{N_1,\cdots,N_d\}$.
The following theorem comes from \cite{WangJinXie}.
\begin{theorem}\cite[Theorem 3]{WangJinXie}\label{theorem_Gauss}
Assume that the function $F(x)$ has the expansion (\ref{def_F(x)}), 
where the coefficient $A_\alpha(x)$ has the form (\ref{def_A_alpha}). 
We use Gauss quadrature points and corresponding weights
to the integration of $F(x)$ on the $d$-dimensional tensor product domain $\Omega$.
If the function $\Psi(x)$ involved in the function $F(x)$ has TNN form (\ref{def_TNN_normed}),
the efficient quadrature scheme (\ref{eq_I_tensor_form}) has $2\underline N$-th order accuracy.
Let $T_1$ denote the computational complexity
for the $1$-dimensional function evaluation operations. 
The computational complexity can be bounded by
$\mathcal O\big(dqT_1k^2p^k\big((d+m)^m+k\big)^kN\big)$, 
which is the polynomial scale of the dimension $d$.
\end{theorem}
\begin{remark}
Other types of one-dimensional quadrature schemes can also be adopted to do
the $d$-dimensional integration.
In numerical examples, we decompose each $\Omega_i$ into subintervals with 
mesh size $h$ and choose $N_i$ one-dimensional Gauss points in each subinterval. 
Then the deduced $d$-dimensional quadrature scheme has accuracy 
$\mathcal O(h^{2\underline N}/(2\underline N)!)$,
where the included constant depends on the smoothness of $F(x)$.
\end{remark}

\begin{remark}
In this paper, we assume $\Omega=\Omega_1\times\cdots\times\Omega_d$ and each $\Omega_i$
is a bounded interval $(a_i,b_i)$.
It is worth mentioning that a similar quadrature scheme can also be given for unbounded $\Omega_i$.
In \cite{WangXie}, in addition to Legendre-Gauss quadrature scheme for bounded domain,
we discuss Hermite-Gauss quadrature scheme for the whole line $\Omega_i=(-\infty,+\infty)$
and Laguerre-Gauss quadrature scheme for the half line $\Omega_i=(a_i,+\infty)$.
The computational complexity of these integrations is also polynomial order of $d$.
\end{remark}

\section{TNN based machine learning method}\label{Section_TNN_ML}
%\subsection{Implementation of TNN}
In this section, we introduce the TNN-based
machine learning method to solve the elliptic
diffusion problem with stochastic diffusion coefficient.
We assumed the physical domain $D=D_1\times\cdots\times D_d$ with $D_i=[a_i, b_i]$
for $i=1, \cdots, d$. It is noted that the tensor product structure plays
an important role in reducing the dependence on dimensions in numerical integration.
We will find that the high accuracy of the high-dimensional integrations of
TNN functions leads to the corresponding machine learning method has high accuracy
for solving the high-dimensional parametric elliptic equations.

For designing the TNN based machine learning method, we build the following TNN function
\begin{eqnarray}\label{tab:TNN_train}
\Psi(y,x)=\sum\limits_{k=1}^p c_k\left(\prod_{j=1}^M\widehat\phi_{j,k}(y_j)\prod_{i=1}^d\widehat\psi_{i,k}(x_i)\right),
\end{eqnarray}
where $y=[y_1, \cdots, y_M]^\top$ and $x=[x_1, \cdots, x_d]^\top$.

In order to treat the boundary condition, following \cite{WangJinXie},
for $i=1,\cdots, d$, the $i$-th subnetwork $\psi_i(x_i;\theta_i)$ is defined as follows:
\begin{eqnarray*}\label{bd_subnetwork}
\psi_i(x_i)&:=&(x_i-a_i)(b_i-x_i)\widehat\psi_i(x_i)\nonumber\\
&=&\big((x_i-a_i)(b_i-x_i)\widehat\psi_{i,1}(x_i),\cdots,(x_i-a_i)(b_i-x_i)
\widehat\psi_{i,p}(x_i)\big)^\top,
\end{eqnarray*}
where $\widehat\psi_i(x_i;\theta_i)$ is an FNN from $\mathbb R$ to $\mathbb R^p$ with sufficiently
smooth activation functions, such that $\Psi(y,x)\in L_\rho^2(H_0^1)$.

The trial function set $V$ is modeled by the TNN structure $\Psi(y,x)$
where the parameters take all the possible values and it is obvious
that $V\subset L_\rho^2(H_0^1)$ by selecting the appropriate activation function.

Based on the Ritz method, the solution $\Psi^*(y,x)$ of the following optimization
problem is the approximation to the problem (\ref{Weak_Problem})
\begin{eqnarray}\label{ML_Weak}
&&\Psi^*(y,x)\nonumber\\
&&= \arg\min\limits_{\Psi(y,x)\in V}
\int_\Gamma \int_D\left[\frac{1}{2}a_M(y,x)\left|\nabla_x\Psi(y,x)\right|^2
- f(y,x)\Psi(y,x)\right]\rho(y)dxdy.
\end{eqnarray}
We can also use the strong form of the stochastic partial differential equation
to build the optimization problem to produce the TNN approximation.
The corresponding optimization problem can be defined as follows
\begin{eqnarray}\label{ML_Strong}
\Psi^*(y,x)=\arg\min\limits_{\Psi(y,x)\in V}
\int_\Gamma \int_D\left|{\rm div}\left(a_M(y,x)\nabla_x\Psi(y,x)\right)
+f(y,x)\right|^2\rho(y)dxdy.
\end{eqnarray}

In practical implementation, the method (\ref{ML_Strong}) always
has better accuracy than the one (\ref{ML_Weak}).
But the required memory for the method (\ref{ML_Strong})
is more than that for the one (\ref{ML_Weak}) since we need to compute the second order
derivatives of the trial function $\Psi(y,x)$.

In this paper, the gradient descent (GD) method is adopted to minimize the loss function
$L(\theta)$. The GD scheme can be described as follows:
\begin{equation}\label{GD_Step}
\theta^{(k+1)}=\theta^{(k)}-\eta\nabla L(\theta^{(k)}),
\end{equation}
where $\theta^{(k)}$ denotes the parameters after the $k$-th GD step, $\eta$ is
the learning rate (step size). In (\ref{GD_Step}), the loss function is defined as
\begin{eqnarray}\label{Loss_Weak}
L:=\int_\Gamma \int_D\left[\frac{1}{2}a_M(y,x)\left|\nabla_x\Psi(y,x)\right|^2
- f(y,x)\Psi(y,x)\right]\rho(y)dxdy,
\end{eqnarray}
for the method (\ref{ML_Weak}) and
\begin{eqnarray}\label{Loss_Strong}
L:=\int_\Gamma \int_D\left|{\rm div}\left(a_M(y,x)\nabla_x\Psi(y,x)\right)
+f(y,x)\right|^2\rho(y)dxdy,
\end{eqnarray}
for the method (\ref{ML_Strong}). In practical learning process,
we use Adam optimizer \cite{KingmaAdam} to choose the learning rate
adaptively and get the optimal solution $\Psi^*(y,x)$.

\section{Numerical examples}\label{Section_Numerical_Examples}
In this section, we provide three numerical examples of solving the elliptic diffusion
problem with stochastic diffusion coefficients by TNN-based machine learning method.
Since the KL eigenpairs can be computed by some classical numerical methods,
such as standard finite element method and spectral method,
here we only consider analytically known KL expansions and their tensorized versions.

Generally speaking, ones are interested in deterministic statistics of $u(\omega, x)$
rather than in the random solution itself. Therefore, in the following numerical examples,
we will show the results regarding the deterministic statistics of $u(\omega, x)$.

To demonstrate the convergence behavior and accuracy of approximations using TNN,
we define the $L^2$ projection operator 
$\mathcal{P}: L^2_\rho(H_0^1) \rightarrow {\rm span}\{\Psi(y,x;\theta^*)\}$ as follows:
\begin{equation*}
\begin{array}{cc}
\langle\mathcal{P}u,v\rangle_{L^2}=\langle u,v\rangle_{L^2}:= \displaystyle\int_{\Gamma}\displaystyle\int_{D}u(y,x)v(y,x)\rho(y)dxdy,\enspace\forall
v\in {\rm span}\left\{\Psi(y,x;\theta^*)\right\},
\end{array}
\end{equation*}
 for $u\in L^2_\rho(H_0^1)$.

And we define the $H^1$ projection operator
$\mathcal{Q}: L^2_\rho(H_0^1) \rightarrow {\rm span}\{\Psi(y,x;\theta^*)\}$ as follows:
\begin{equation*}
\begin{array}{cc}
\langle\mathcal{Q}u,v\rangle_{H^1}=\langle u,v\rangle_{H^1}:=\displaystyle\int_{\Gamma}\displaystyle\int_{D}\nabla_x u(y,x)\cdot\nabla_x v(y,x)\rho(y)dxdy,\enspace\forall v\in
{\rm span}\{\Psi(y,x;\theta^*)\},
\end{array}
\end{equation*}
for $u\in L^2_\rho(H_0^1)$.

Then we define the following relative errors for the approximated function $\Psi(y,x;\theta^*)$
\begin{equation}
\begin{array}{c}
e_{L^2}:=\displaystyle\frac{\Vert u-\mathcal{P}u\Vert_{L^2_\rho(\Gamma;L^2(D))}} 
{\Vert f\Vert_{L^2(D)}},
\enspace\enspace e_{H^1}:=
\displaystyle\frac{\|u-\mathcal{Q}u\|_{L^2_\rho(\Gamma;H^1(D))}}{\vert f\vert_{H^1(D)}}.
\end{array}
\end{equation}

\subsection{Zero-order coefficients}
In the first example, we consider the case of spatially homogeneous stochastic coefficients,
\begin{eqnarray}
a_M(y,x)=1+\sum\limits_{m=1}^M a_my_m,
\end{eqnarray}
for the truncated sequence of (spatially homogeneous) coefficients
$a_m=(1+m)^{-\alpha}(m=1,\cdots,M)$ with algebra decay rates $\alpha=2$.
Here, we set $D=(0,1)$.
The variables $y_m$ $(m=1,2\cdots,M)$ are independently and identically
distributed and follow the uniform distribution on $[-1,1]$.
Here, we wish to construct an equation that has the analytical solution
$u(y,x)=\sin(\pi x)\prod\limits_{m=1}^M \sin\left(\frac{\pi}{2} y_m\right)$.
For this aim,  the right hand side load is set to be
\begin{eqnarray}
f(y,x) = \left(1+\sum\limits_{i=1}^M a_i y_i\right)\cdot\pi^2 \sin(\pi x)
\prod\limits_{m=1}^M \sin\left(\frac{\pi}{2}y_m\right).
\end{eqnarray}

In this example, we choose $M=10$, $20$, $50$ and $100$,
the computational probability domain is set to be
$\Gamma=\Gamma_1\times\cdots\times\Gamma_M$, where $\Gamma_i=[-1,1]$ for $i=1,\cdots, M$.
In order to do the high-dimensional integration for the TNN functions,
we decompose the physical domain $D$ and the probability
space $\Gamma_j$\  $(j=1,\cdots,M)$ into $200$ subintervals and
choose $16$ Gauss points on each subinterval respectively.

Each subnetwork in TNN is set to the FNN with 3 hidden layers, 
$100$ nodes in each layer is adopted to build the subnetwork for the
TNN function. Here we select $\sin(x)$ as the activation function and choose $p=50$.

For $M=10,20$, we use the strong form to build the loss function (\ref{Loss_Strong}) and
solve the optimization problem (\ref{ML_Strong}) to obtain the approximation of the 
problem (\ref{Weak_Problem}).  The Adam optimizer with
learning rate 0.0005 is adopted for the first 100,000 steps.
Then the LBFGS optimizer with learning rate 0.5 is carried out for the next 10,000 steps.
For the case of $M=50, 100$,  we use the weak form to build the loss function (\ref{Loss_Weak})
and solve the corresponding optimization problem (\ref{ML_Weak}) to produce  the approximation 
of the problem (\ref{Weak_Problem}). The Adam optimizer with
learning rate 0.0001 is adopted for the first 95,000 steps.
Then the LBFGS with learning rate 0.1 is carried out for the next 5,000 steps.

\begin{figure}[htb!]
\centering
\includegraphics[width=3.5cm,height=3.5cm]{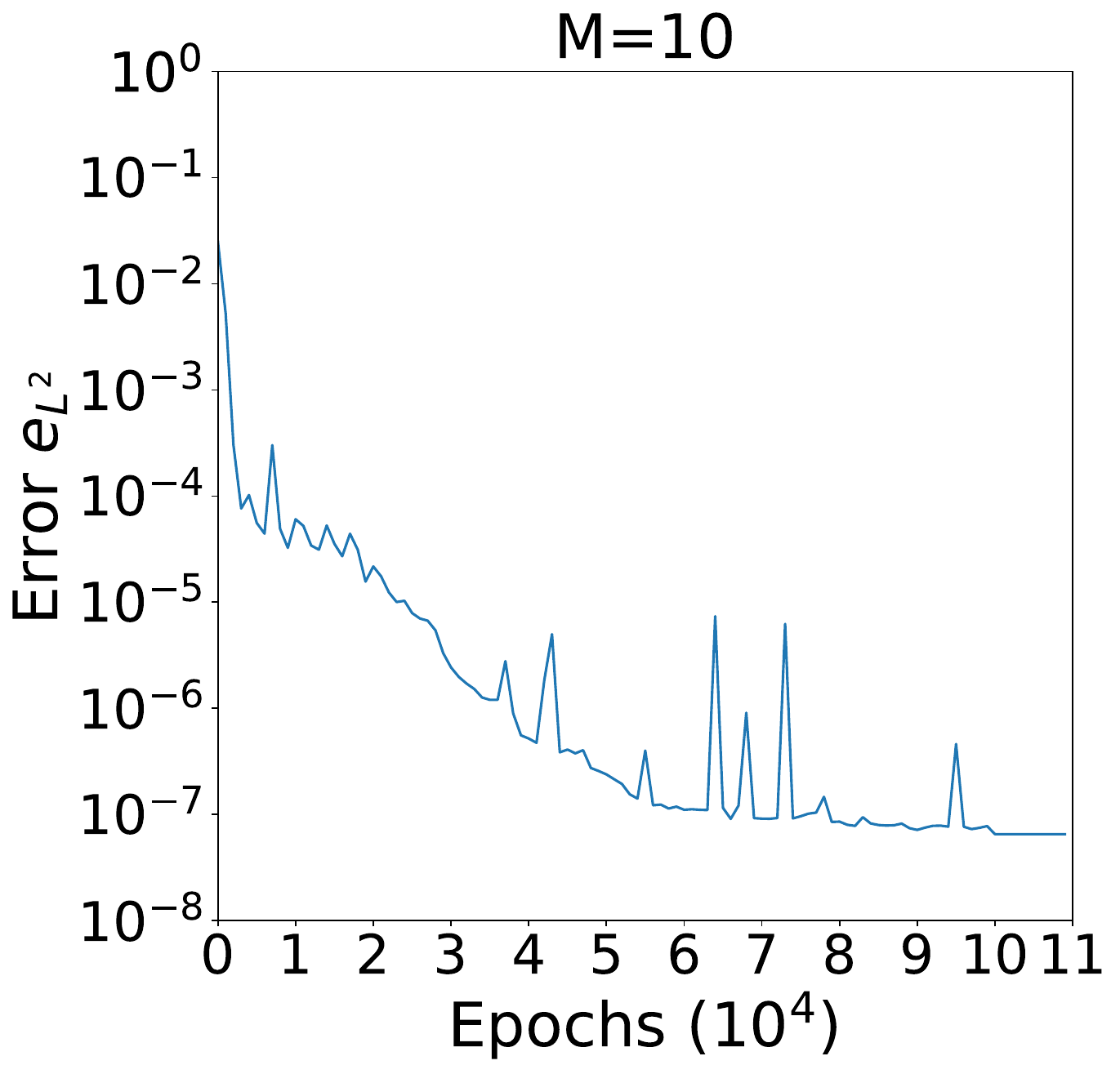}
\includegraphics[width=3.5cm,height=3.5cm]{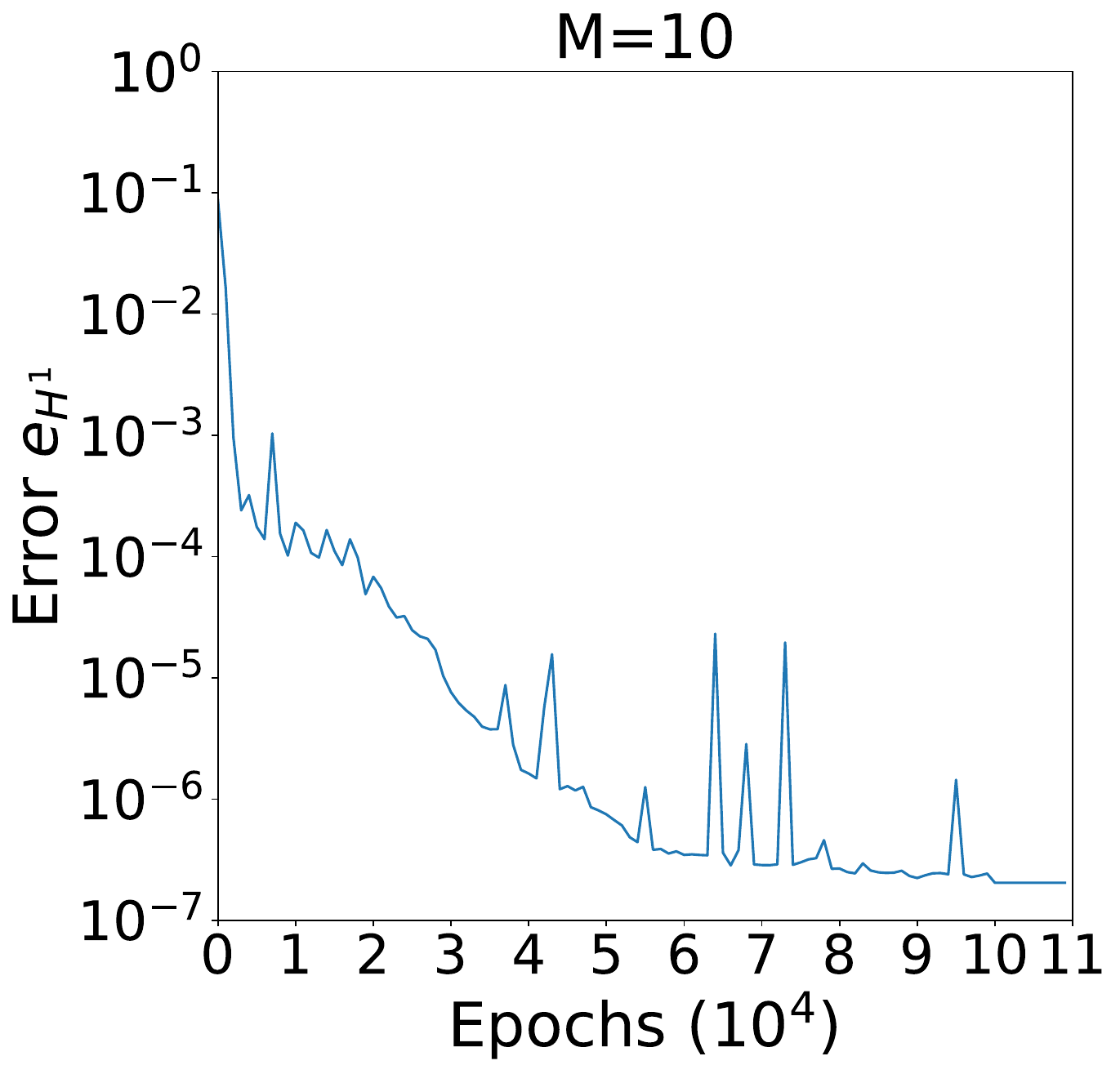}
\includegraphics[width=3.5cm,height=3.5cm]{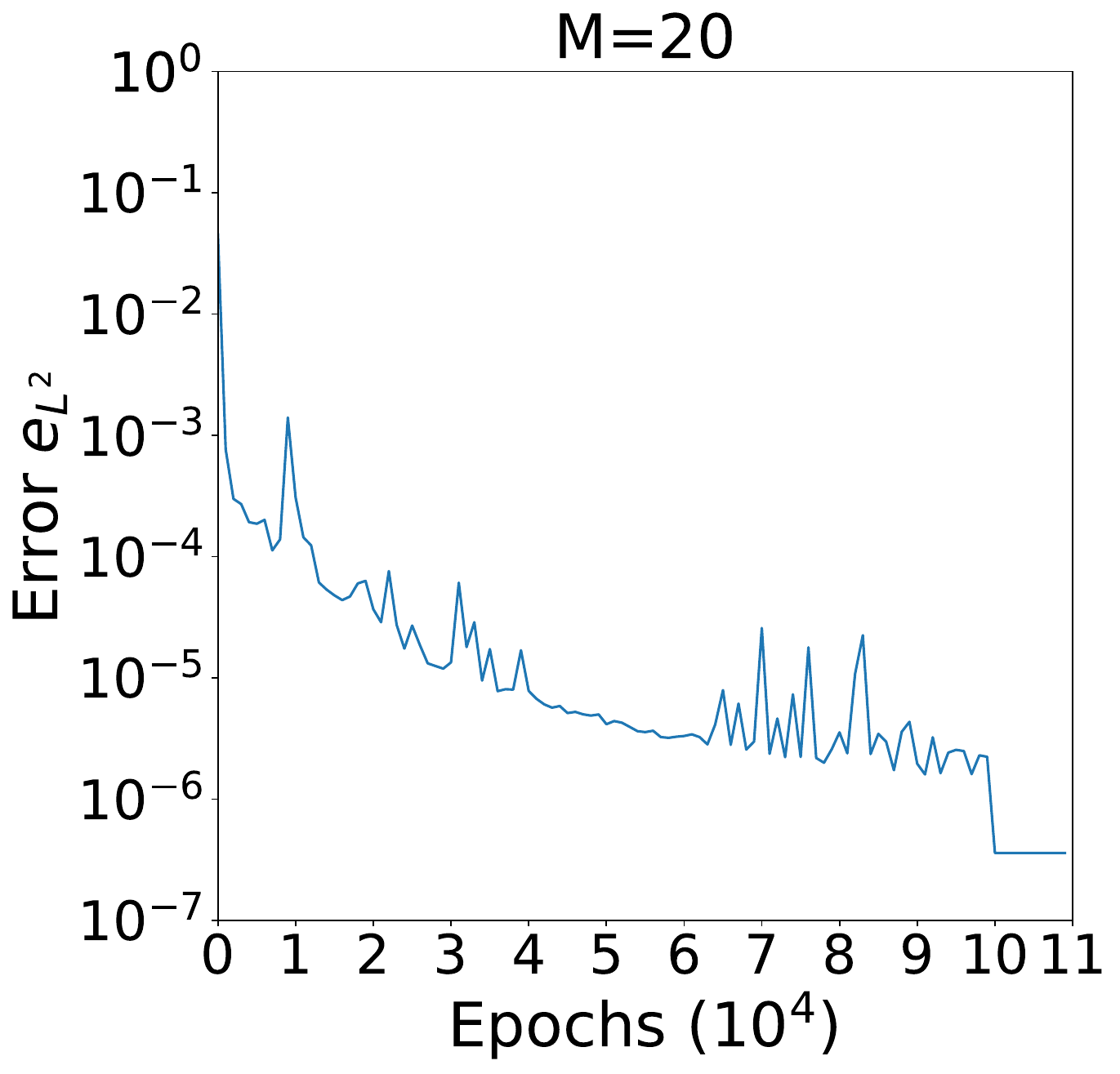}
\includegraphics[width=3.5cm,height=3.5cm]{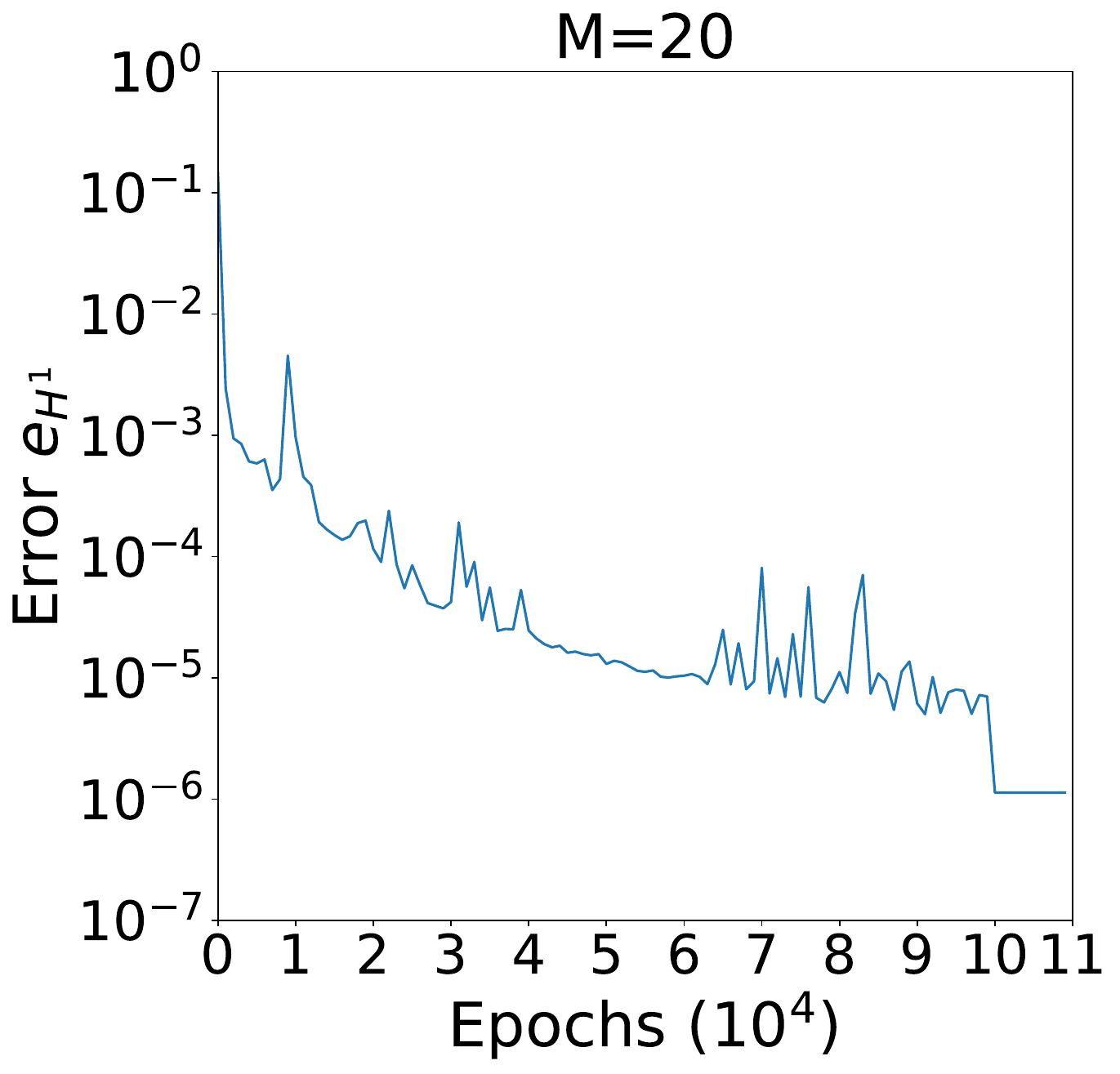}\\
\includegraphics[width=3.5cm,height=3.5cm]{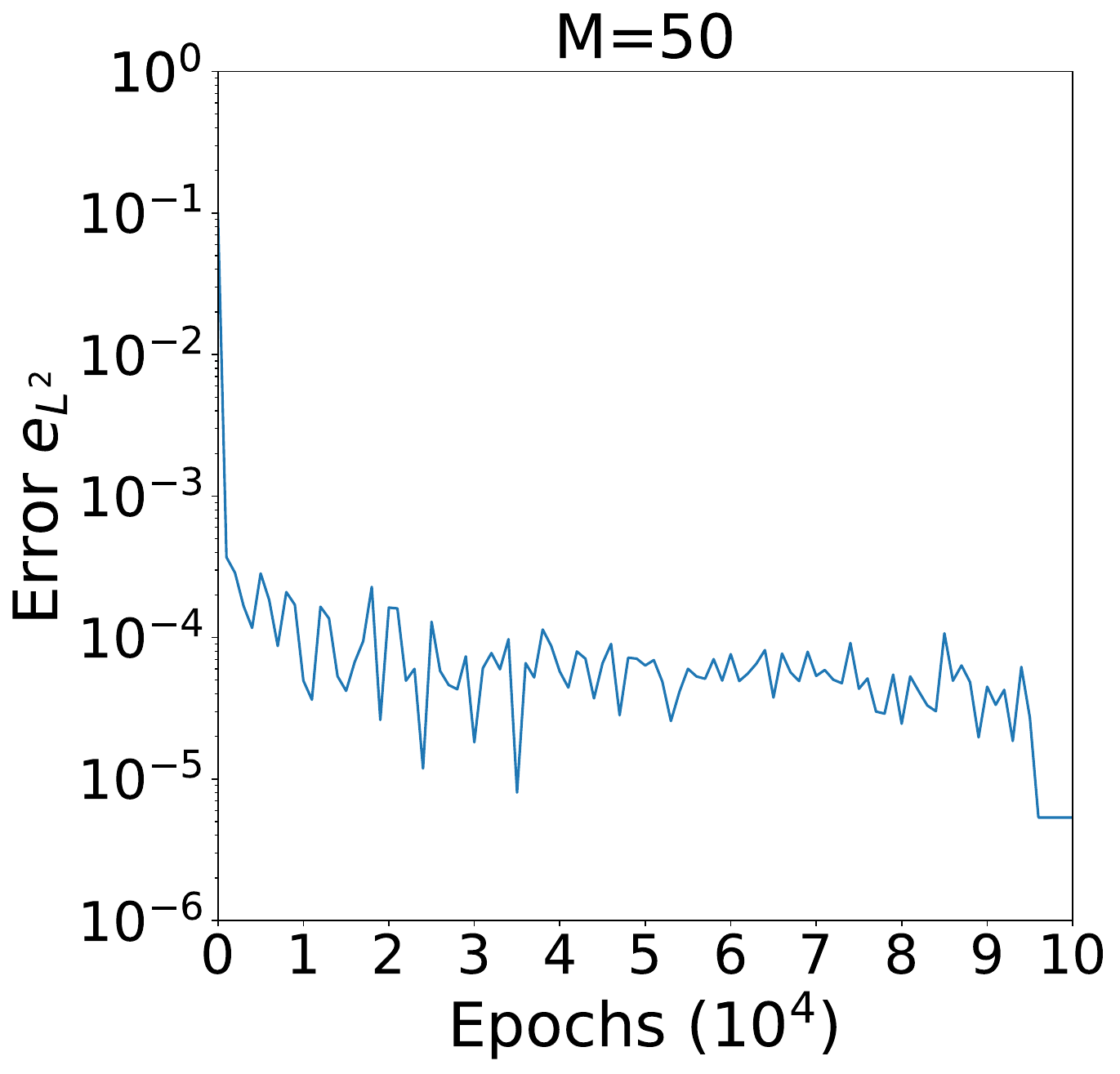}
\includegraphics[width=3.5cm,height=3.5cm]{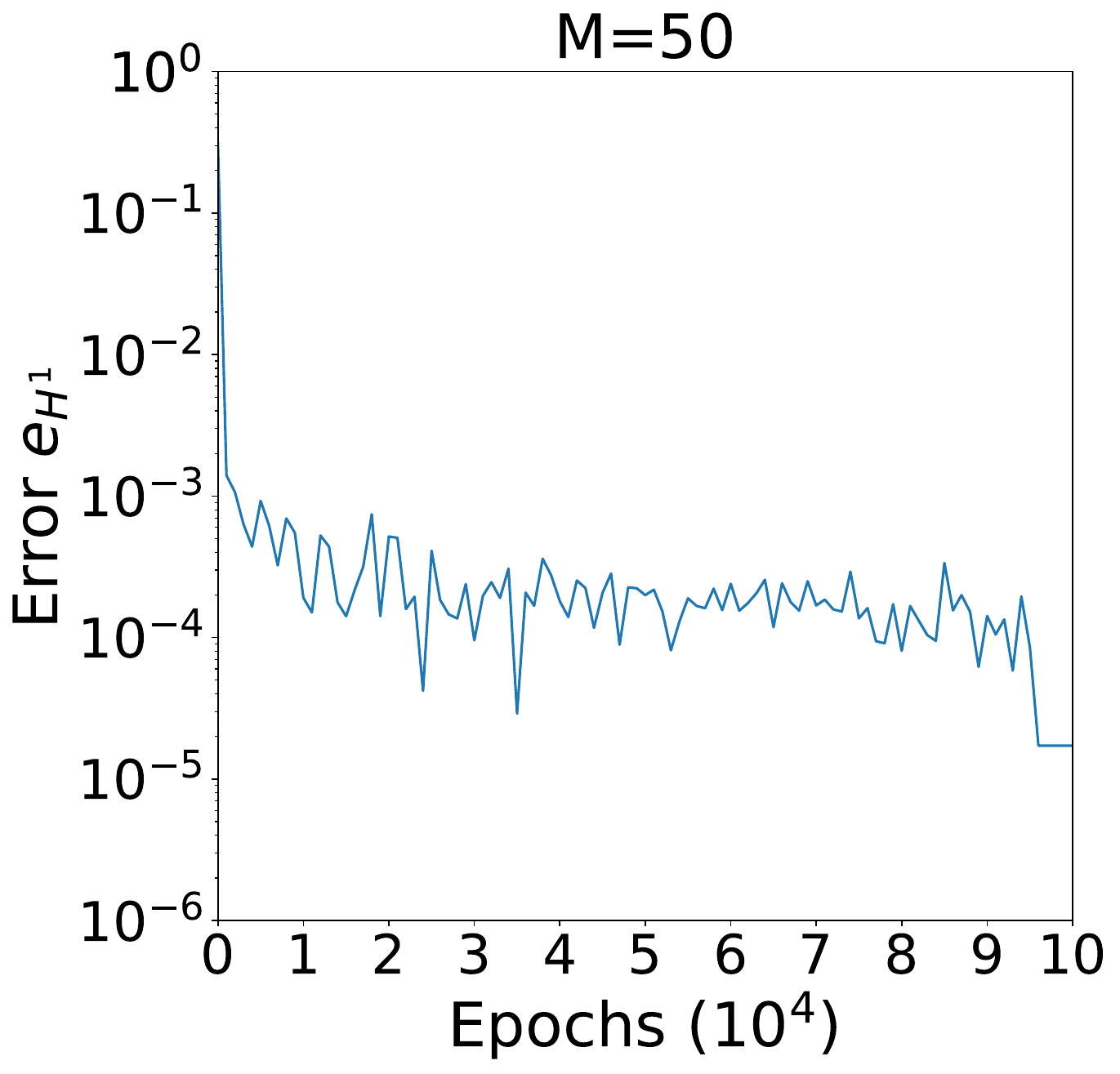}
\includegraphics[width=3.5cm,height=3.5cm]{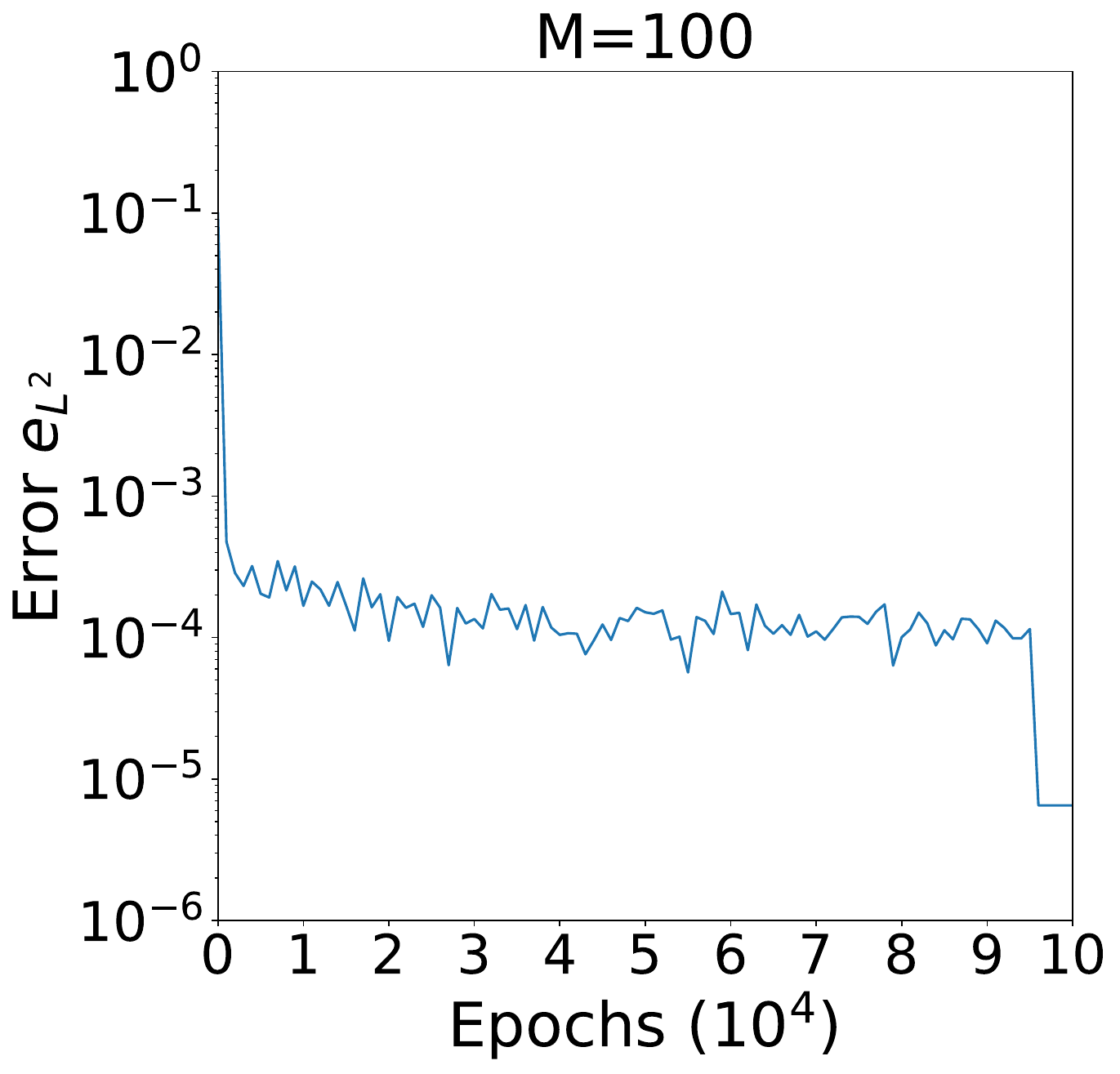}
\includegraphics[width=3.5cm,height=3.5cm]{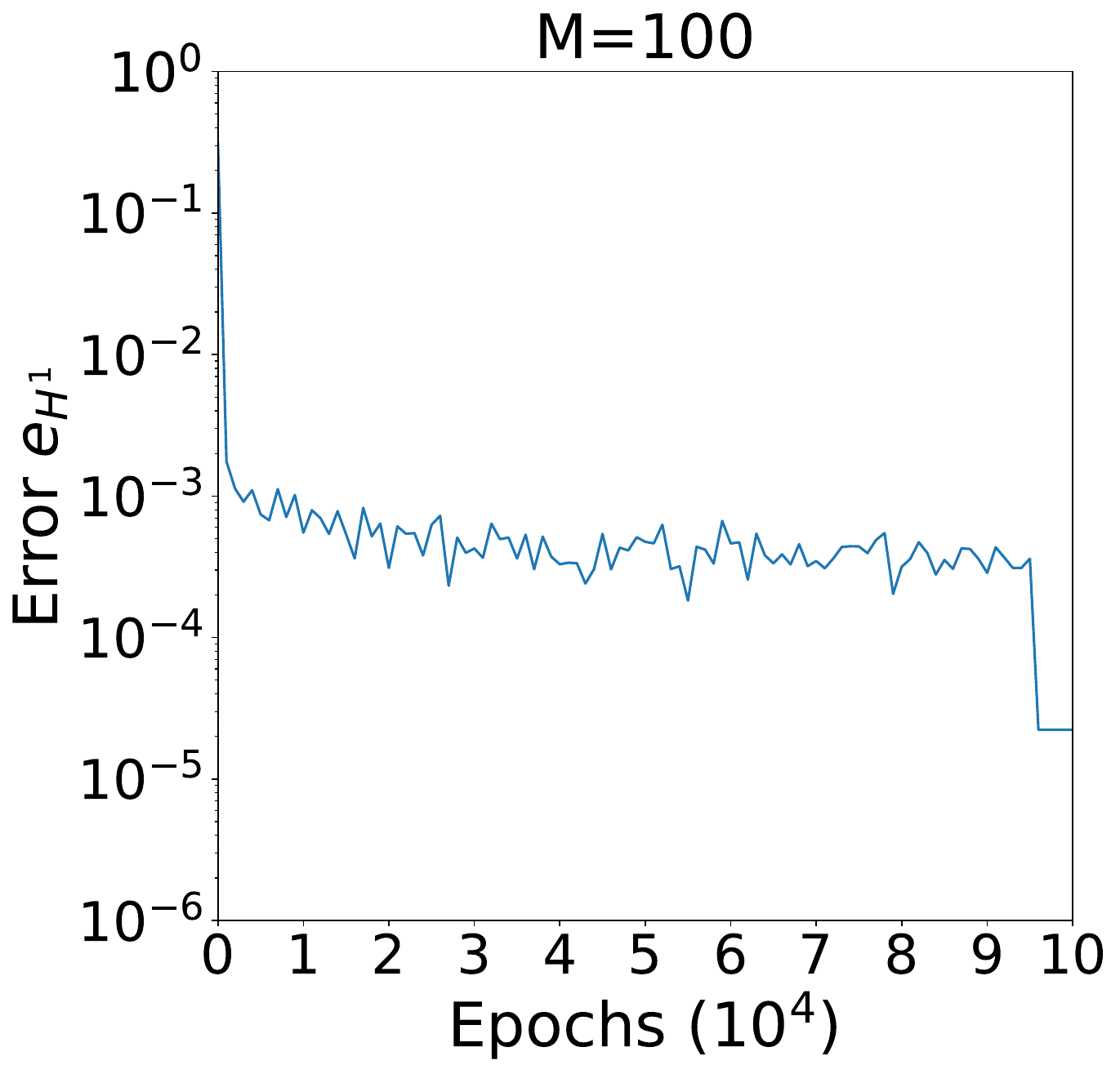}
\caption{The relative errors during training process in Example 1.}\label{fig_errors_ex1}
\end{figure}

The corresponding final errors are listed in the Table \ref{table_errors_ex1}, where we can
find the TNN based machine learning method can solve the high-dimensional stochastic
partial differential equations with high accuracy.   Figure \ref{fig_errors_ex1} 
shows the relative errors $e_{L^2}$ and $e_{H^1}$ versus the number of epochs.
\begin{table}[htb!]
\caption{Errors of Example 1.}\label{table_errors_ex1}
\begin{center}
\vskip-0.2cm
\begin{tabular}{ccccc}
\hline
$M$&   $e_{L^2}$&   $e_{H^1}$\\
\hline
10&   6.515e-08&   2.047e-07\\
20&   3.606e-07&   1.133e-06\\
50&   5.345e-06&   1.720e-05\\
100&  6.506e-06&   2.226e-05\\
\hline
\end{tabular}
\end{center}
\end{table}
%------------------------------------
%========================================================================================
\subsection{Variable stochastic coefficients}
In the second example, we consider the case of variable stochastic coefficients
\begin{eqnarray}
a_M(y,x)=1+\sum\limits_{m=1}^M a_m(x)y_m
\end{eqnarray}
for the truncated sequence of coefficients $a_m(x)=(1+m)^{-2}\sin(m\pi x)$
$(m=1,\cdots,M)$.
Here, we also set $D=(0,1)$.

Similarly, we assume that the exact solution
$$u(y,x)=\sin(\pi x)\prod\limits_{m=1}^M \sin\left(\frac{\pi y_m}{2}\right)$$
and the load is chosen accordingly
\begin{eqnarray*}
f(y,x)&=&\pi^2\sin(\pi x)\sum_{i=1}^M\sin\left(\frac{\pi y_i}{2}\right) + \sum_{m=1}^M\left(\frac{\pi}{1+m}\right)^2\sin(m\pi x)\sin(\pi x)y_m
\prod_{i=1}^M\sin\left(\frac{\pi y_i}{2}\right)\nonumber\\
&&-\sum_{m=1}^M\frac{m\pi^2}{(1+m)^2}\cos(m\pi x)\cos(\pi x)y_m\prod_{i=1}^M
\sin\left(\frac{\pi y_i}{2}\right).
\end{eqnarray*}

In this example, we investigate the numerical performance for the cases 
$M= 10$, $20$, $50$ and $100$. 
Here, the computational probability domain is also set to be
$\Gamma=\Gamma_1\times\cdots\times\Gamma_M$, where $\Gamma_i=[-1,1]$ for $i=1,\cdots, M$.
In order to do the high-dimensional integration for the TNN functions,
we also decompose the physical domain $D$ and the probability
space $\Gamma_j$\  $(j=1,\cdots,M)$ into $200$ subintervals and
choose $16$ Gauss points on each subinterval respectively.

The FNN with 3 hidden layers, $100$ nodes in each layer is also adopted to 
build the subnetwork for the TNN function. Here we select $\sin(x)$ 
as the activation function and choose $p=50$.

For $M=10, 20$, we also use the strong form to build the loss function (\ref{Loss_Strong}) and
solve the optimization problem (\ref{ML_Strong}) to obtain the approximation 
of the problem (\ref{Weak_Problem}).  The Adam optimizer with
learning rate 0.0005 is adopted for the first 100,000 steps.
Then the LBFGS optimizer with learning rate 0.5 is carried out for the next 10,000 steps.
For the case of $M=50, 100$,  we also use the weak form to build the 
loss function (\ref{Loss_Weak}) and solve the corresponding optimization 
problem (\ref{ML_Weak}) to produce  the approximation of the problem (\ref{Weak_Problem}).
The Adam optimizer with
learning rate 0.0001 is adopted for the first 95,000 steps.
Then the LBFGS with learning rate 0.1 is carried out for the next 5,000 steps.

\begin{figure}[htb!]
\centering
\includegraphics[width=3.5cm,height=3.5cm]{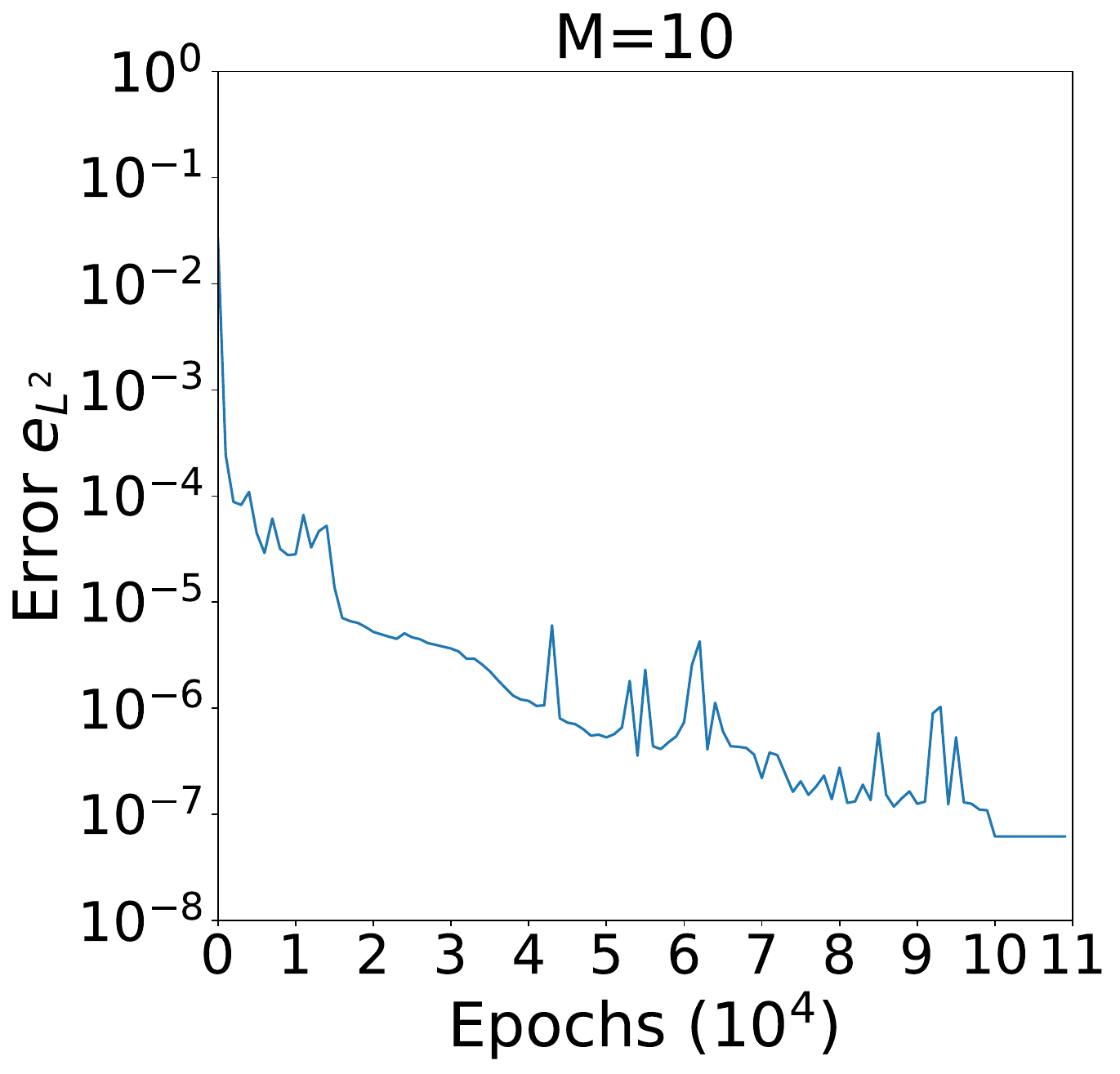}
\includegraphics[width=3.5cm,height=3.5cm]{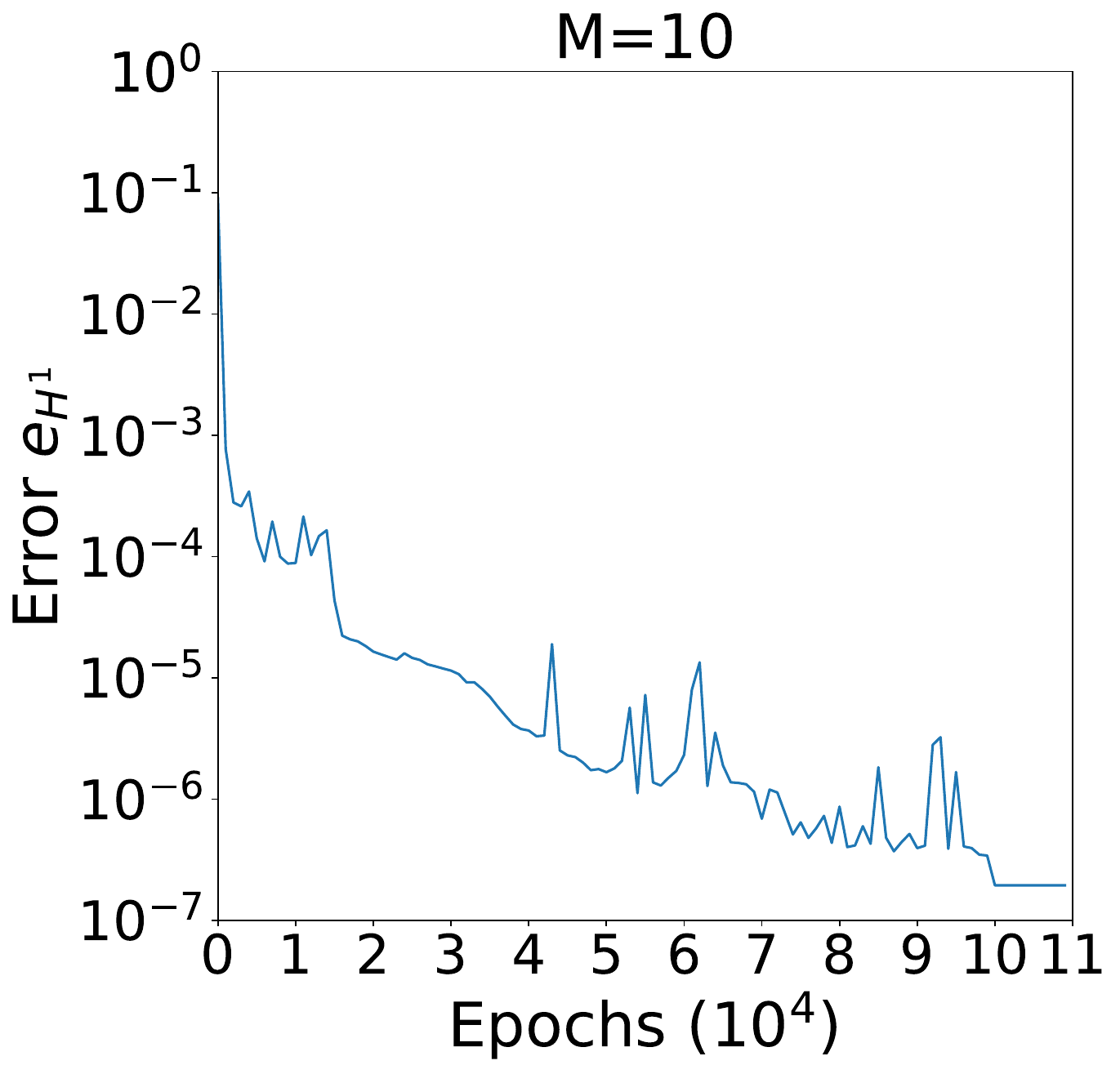}
\includegraphics[width=3.5cm,height=3.5cm]{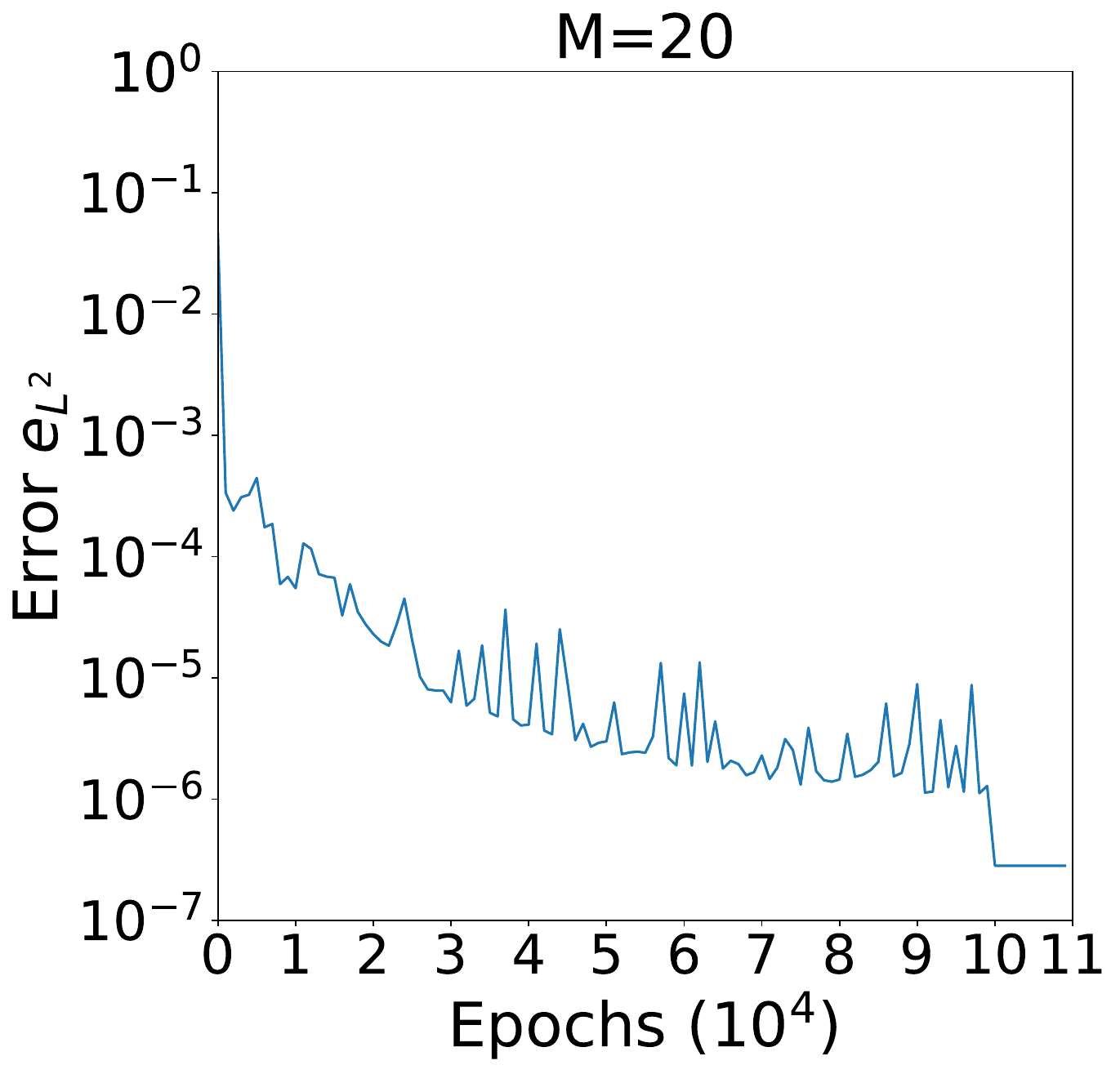}
\includegraphics[width=3.5cm,height=3.5cm]{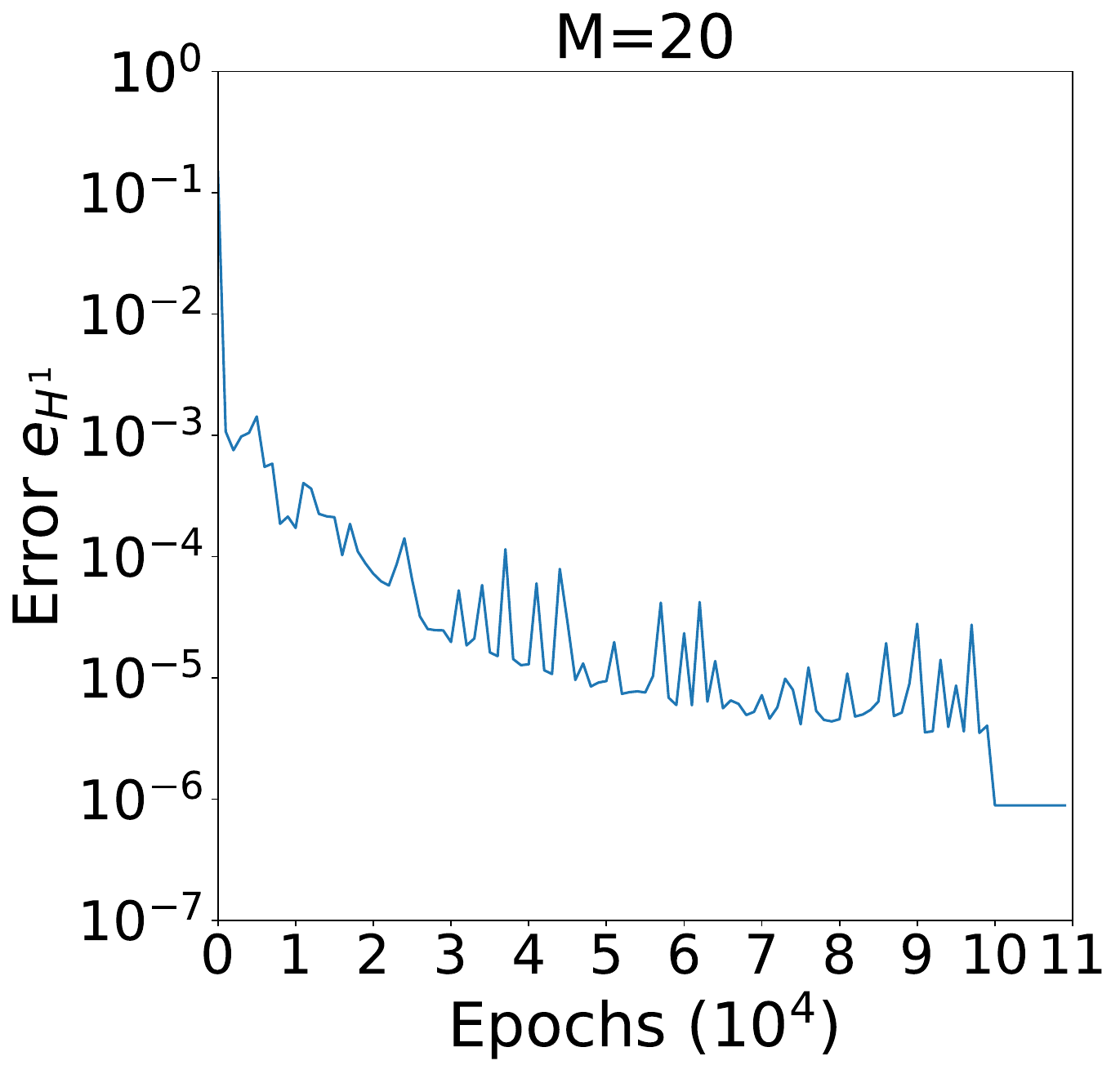}\\
\includegraphics[width=3.5cm,height=3.5cm]{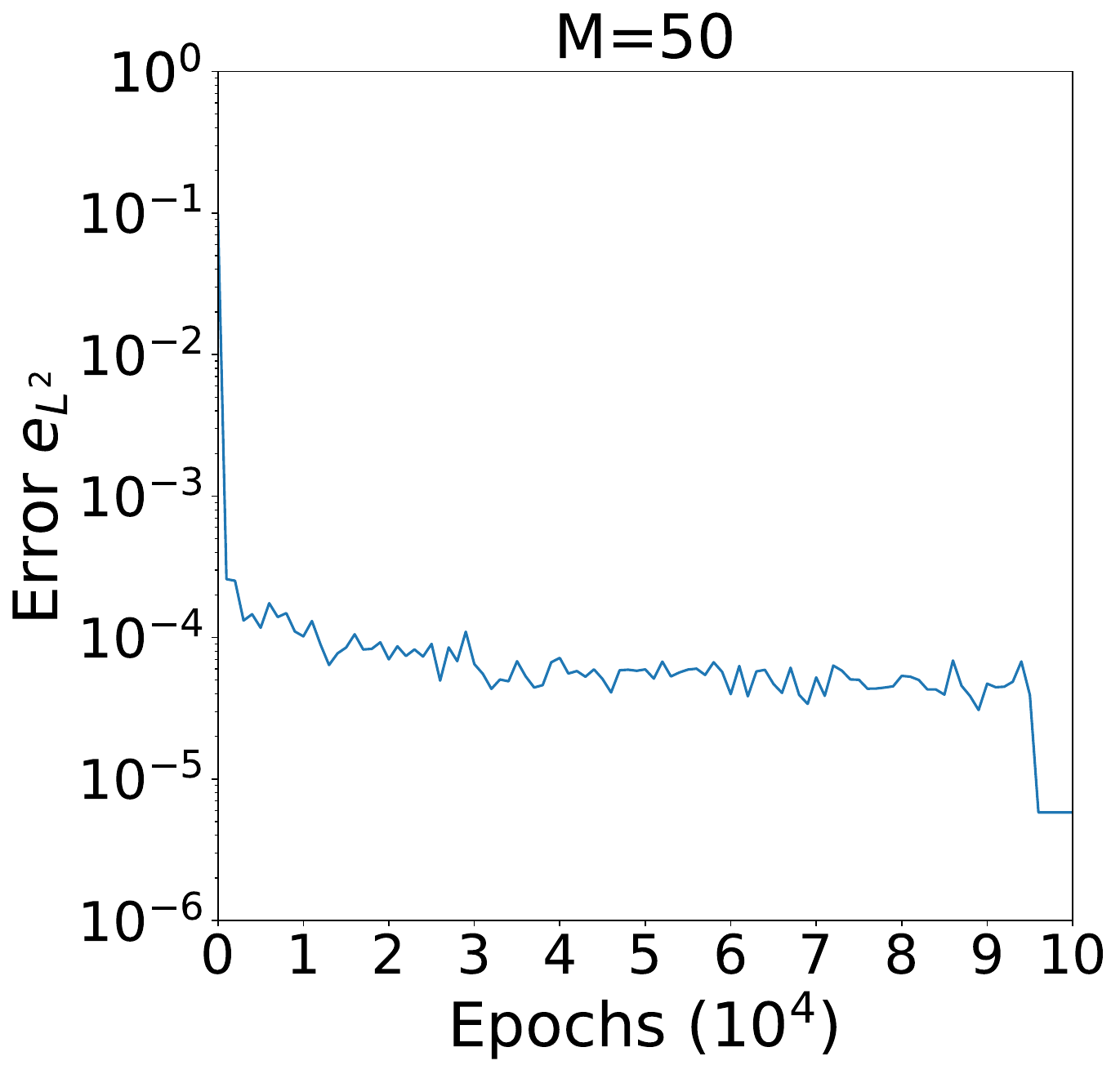}
\includegraphics[width=3.5cm,height=3.5cm]{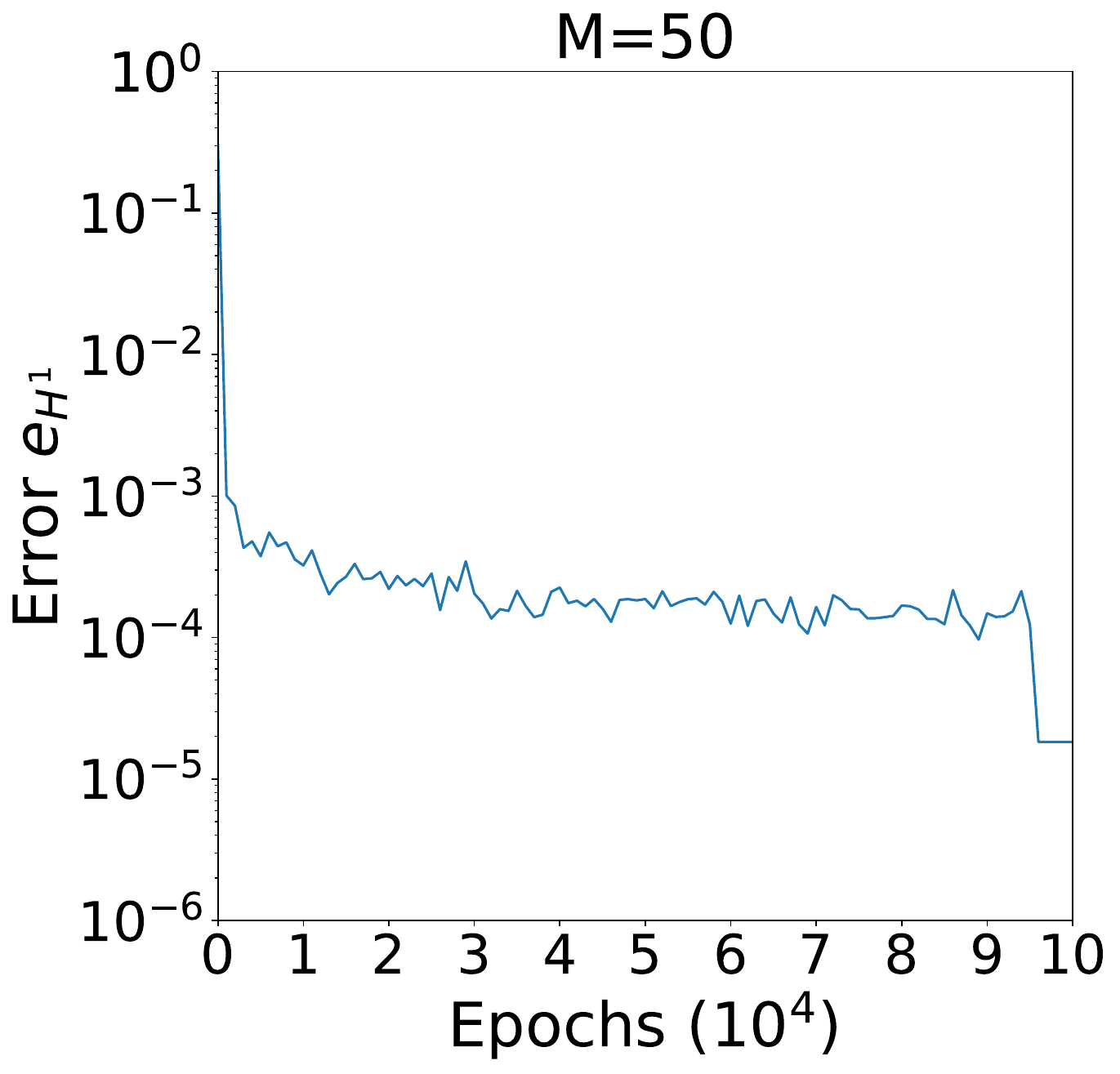}
\includegraphics[width=3.5cm,height=3.5cm]{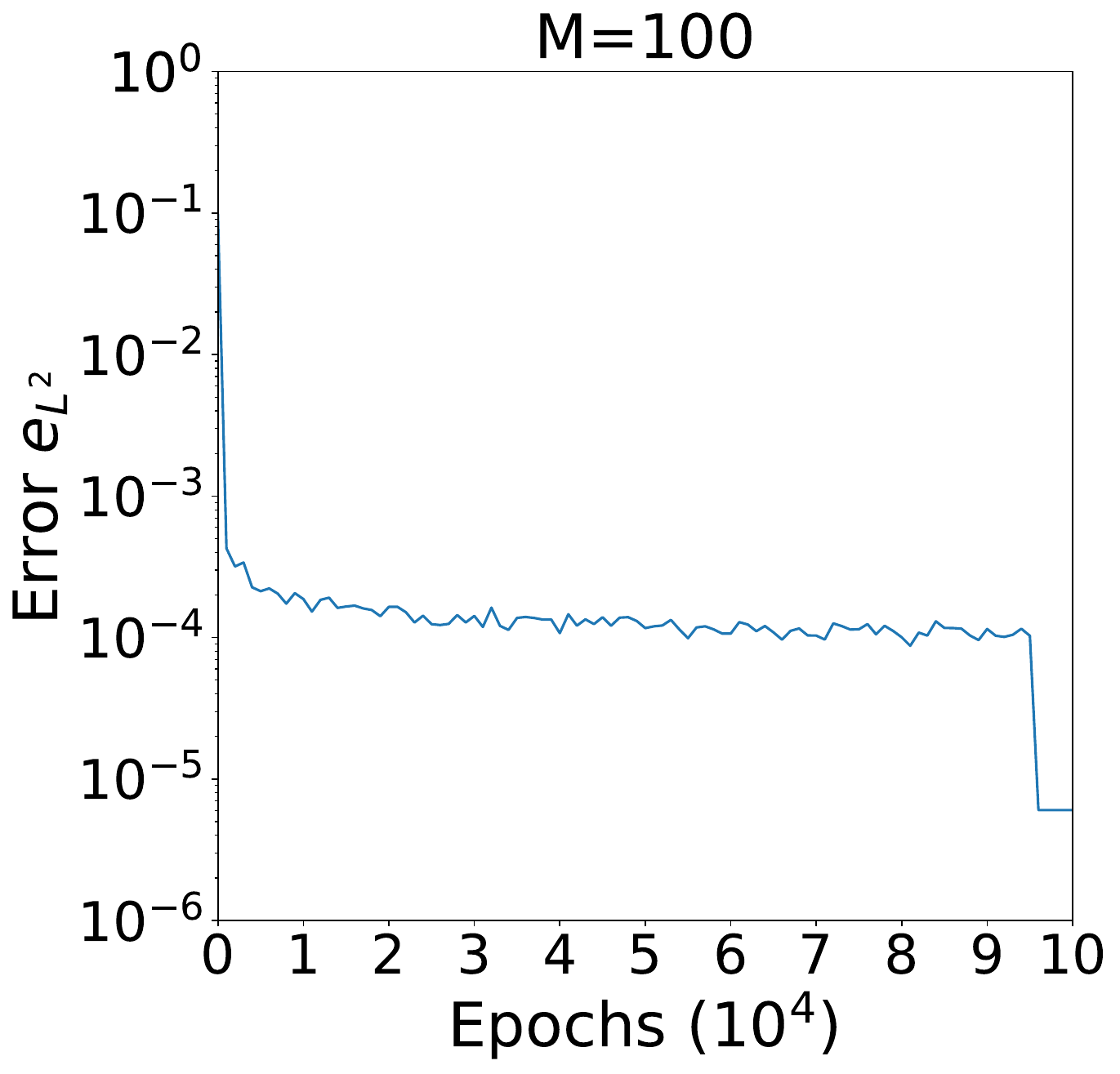}
\includegraphics[width=3.5cm,height=3.5cm]{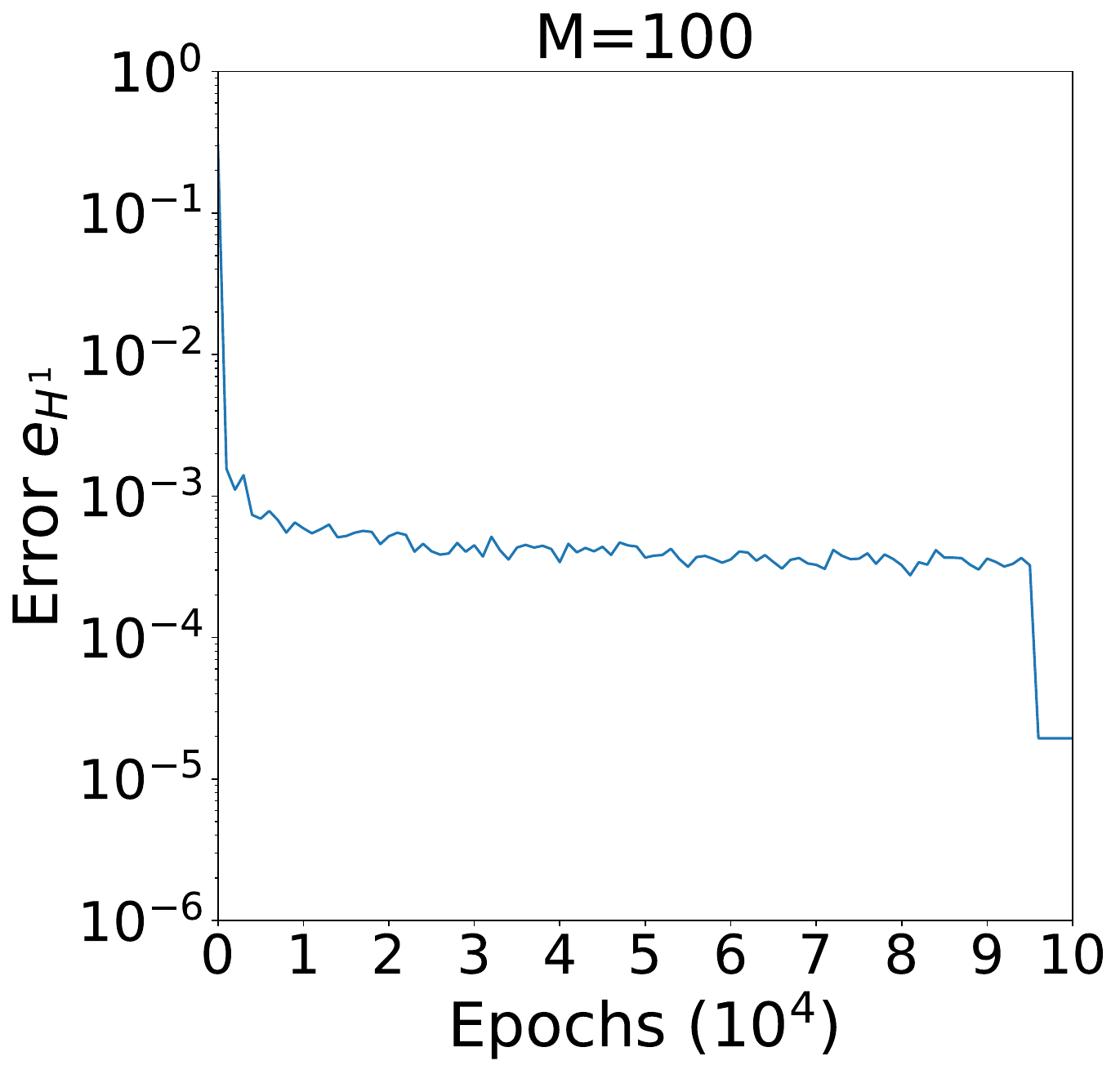}
\caption{The relative errors during training process in Example 2.}\label{fig_errors_ex2}
\end{figure}

The corresponding final errors are listed in the Table \ref{table_errors_ex2}, where we can
find the TNN based machine learning method can solve the high-dimensional stochastic
partial differential equations with high accuracy.   Figure \ref{fig_errors_ex2} 
shows the relative errors $e_{L^2}$ and $e_{H^1}$ versus the number of epochs.
\begin{table}[htb!]
\caption{Errors of Example 2.}\label{table_errors_ex2}
\begin{center}
\begin{tabular}{ccccc}
\hline
$M$ &  $e_{L^2}$ &   $e_{H^1}$\\
\hline
10  &  6.206e-08 &   1.950e-07\\
20  &  2.832e-07 &   8.898e-07\\
50  &  5.803e-06 &   1.830e-05\\
100 &  6.031e-06 &   1.941e-05\\
\hline
\end{tabular}
\end{center}
\end{table}

%========================================================================================
\subsection{Exponential decay coefficients}
In the third example, the variable stochastic coefficients is set to be 
exponential decay coefficient, i.e., 
\begin{eqnarray}
a_M(y,x)=1+\sum\limits_{m=1}^M 0.5\exp(-m)\sin(m\pi x)y_m.
\end{eqnarray} 
Here, we also set $D=(0,1)$. 
Here, in order to make the exact solution is 
$$u(y,x)=\sin(\pi x)\prod\limits_{m=1}^M \sin\left(\frac{\pi y_m}{2}\right),$$ 
we set the load to be 
\begin{eqnarray*}
f(y,x)&=&\pi^2\sin(\pi x)\sum_{i=1}^M\sin\left(\frac{\pi y_i}{2}\right) + 
0.5\pi^2\sum_{m=1}^M\exp(-m)\sin(m\pi x)\sin(\pi x)y_m\prod_{i=1}^M
\sin\left(\frac{\pi y_i}{2}\right)\nonumber\\
&&-0.5\pi^2\sum_{m=1}^Mm\exp(-m)\cos(m\pi x)\cos(\pi x)y_m\prod_{i=1}^M
\sin\left(\frac{\pi y_i}{2}\right).
\end{eqnarray*}

In this example, the numerical experiments are performed to investigate 
the accuracy for the cases $M= 10$, $20$, $50$ and $100$. 
The computational probability domain is also 
$\Gamma=\Gamma_1\times\cdots\times\Gamma_M$ with $\Gamma_i=[-1,1]$ for $i=1,\cdots, M$.
In order to do the high-dimensional integration for the TNN functions, 
the physical domain $D$ and the probability
space $\Gamma_j$\  $(j=1,\cdots,M)$ are decomposed into $200$ subintervals and 
$16$ Gauss points are used on each subinterval.

Here, the FNN with 3 hidden layers, $100$ nodes in each layer is also adopted to 
build the subnetwork for the TNN function. 
Here we select $\sin(x)$ as the activation function and choose $p=50$.

For $M=10, 20$, we use the loss function (\ref{Loss_Strong}) and 
the optimization problem (\ref{ML_Strong}) 
to obtain the approximation of the problem (\ref{Weak_Problem}). 
The Adam optimizer with learning rate 0.0005 is adopted for the first 100,000 steps 
and the LBFGS optimizer with learning rate 0.5 for the sequent 10,000 steps.
For the case of $M=50, 100$,  we use the loss function (\ref{Loss_Weak})
and the corresponding optimization problem (\ref{ML_Weak}) to produce  
the approximation of the problem (\ref{Weak_Problem}). The Adam optimizer with
learning rate 0.0001 is adopted for the first 95,000 steps and 
the LBFGS with learning rate 0.1 for the sequent 5,000 steps.

\begin{figure}[htb!]
\centering
\includegraphics[width=3.5cm,height=3.5cm]{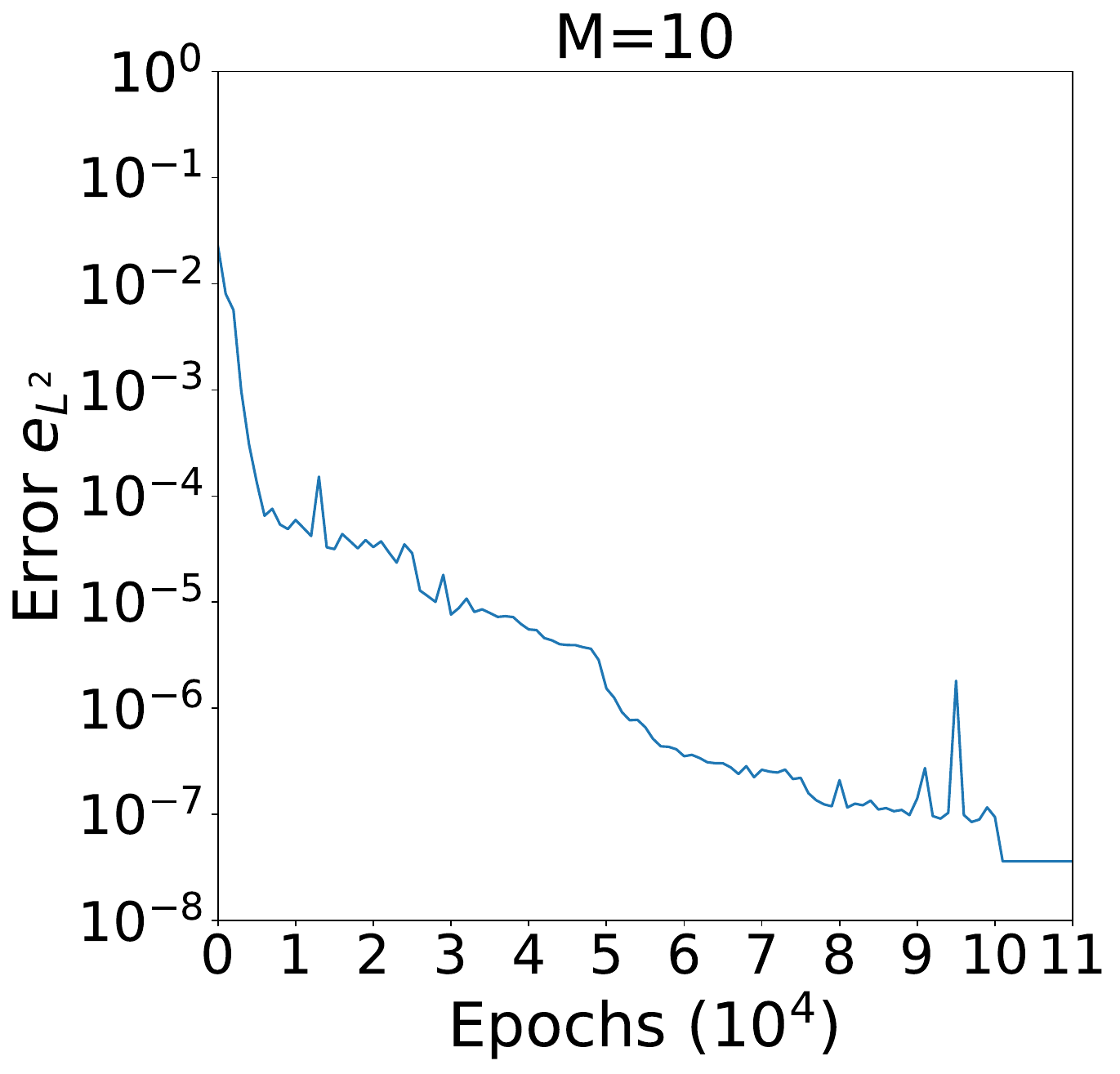}
\includegraphics[width=3.5cm,height=3.5cm]{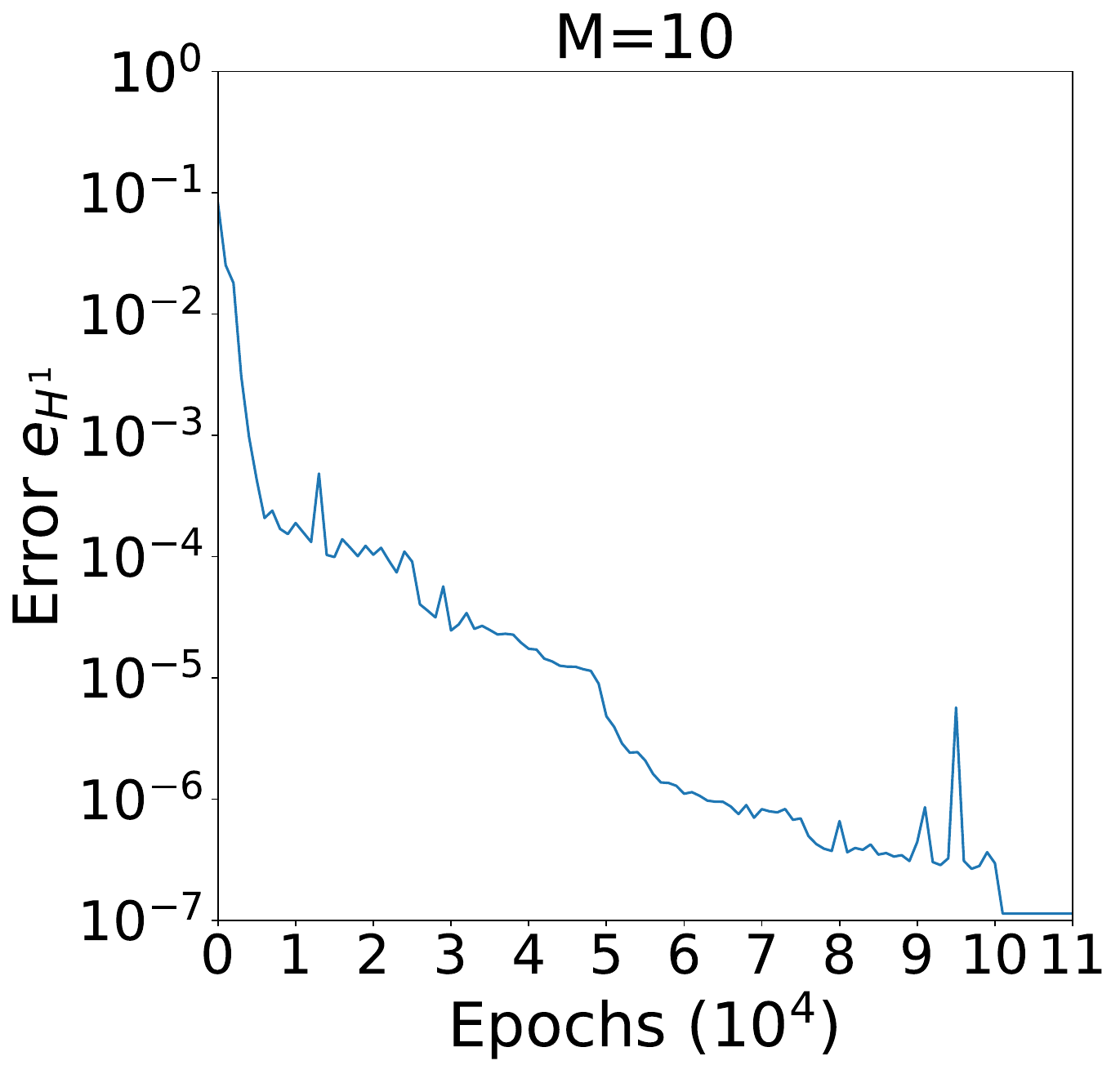}
\includegraphics[width=3.5cm,height=3.5cm]{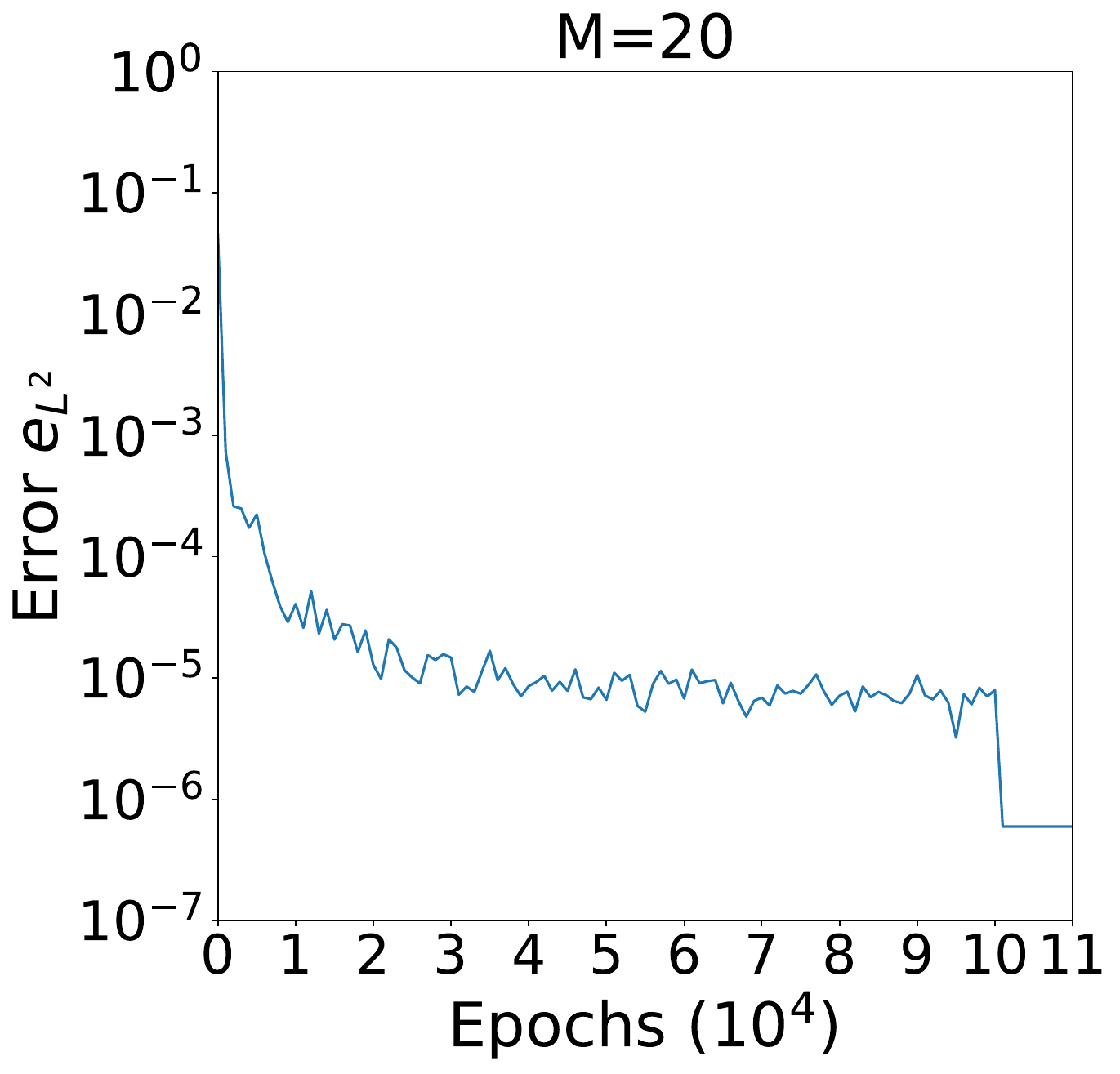}
\includegraphics[width=3.5cm,height=3.5cm]{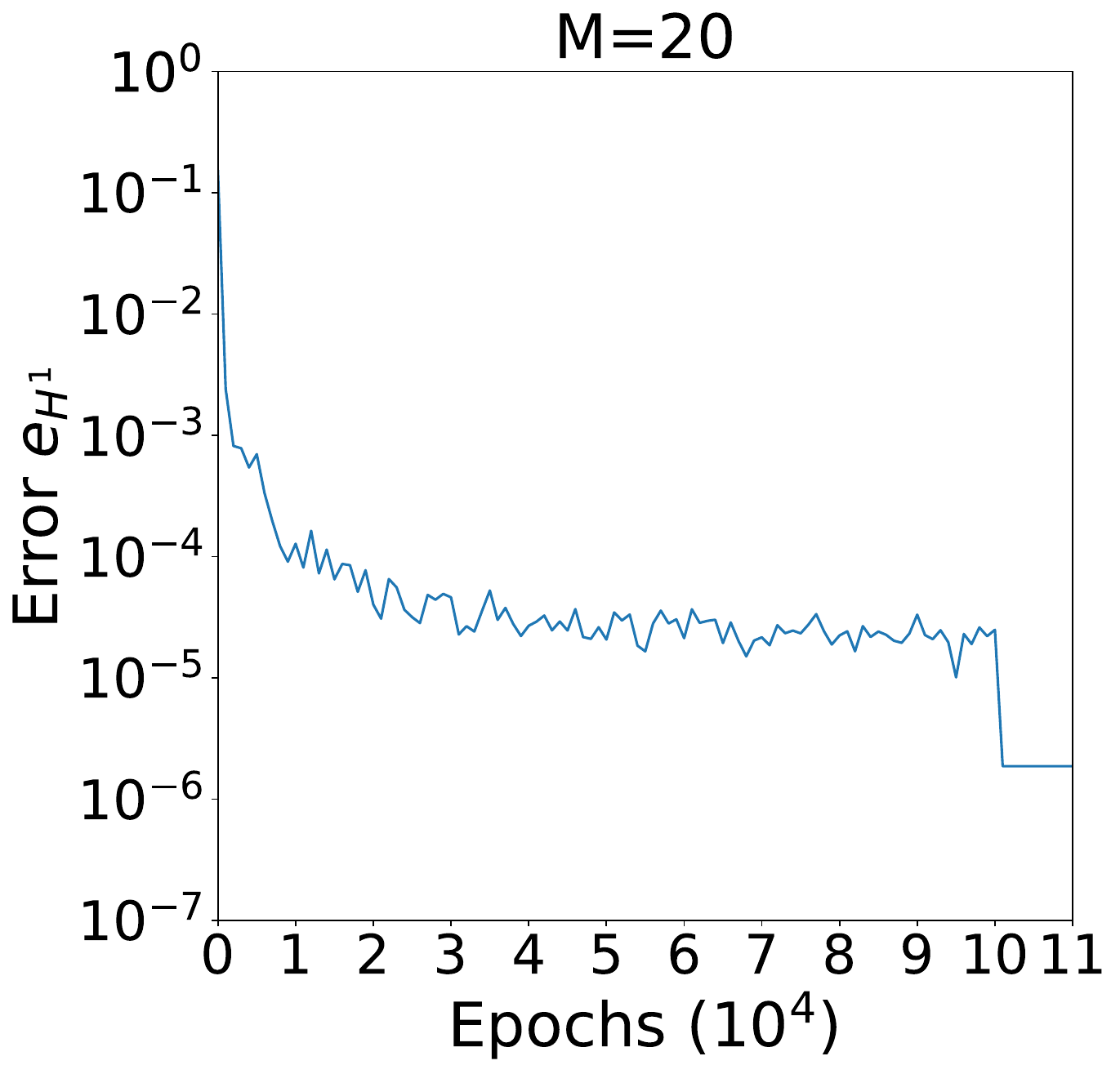}\\
\includegraphics[width=3.5cm,height=3.5cm]{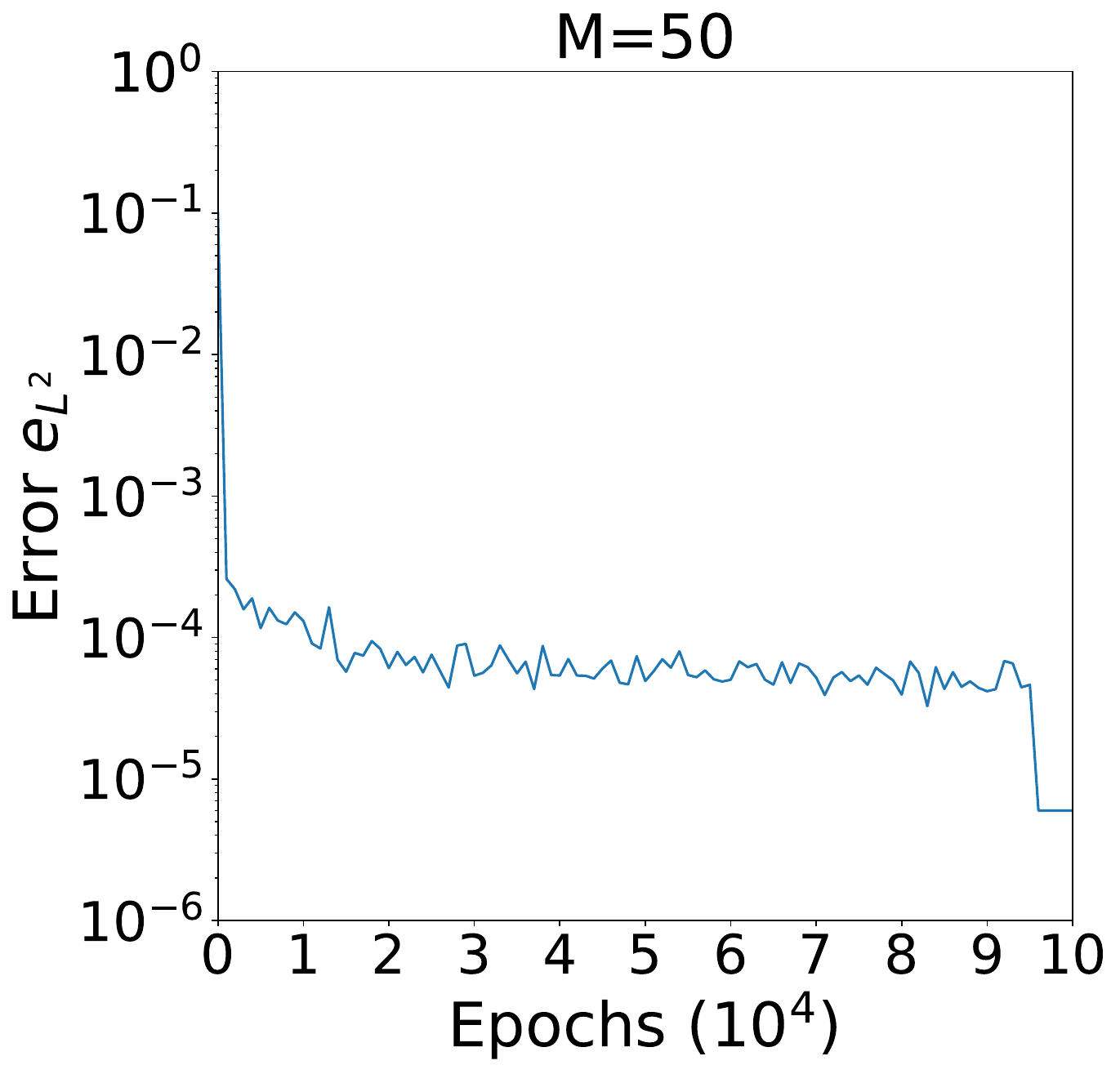}
\includegraphics[width=3.5cm,height=3.5cm]{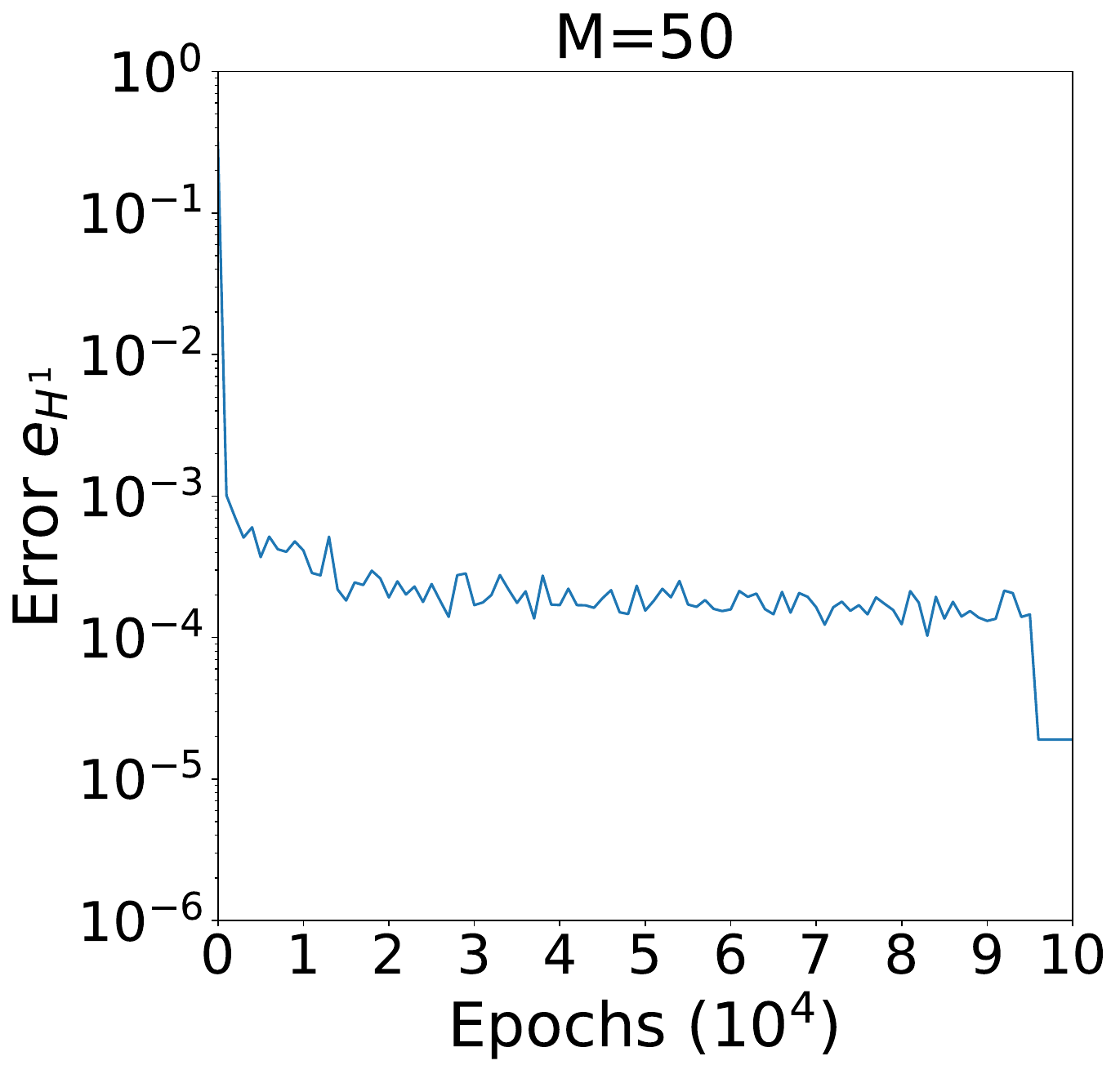}
\includegraphics[width=3.5cm,height=3.5cm]{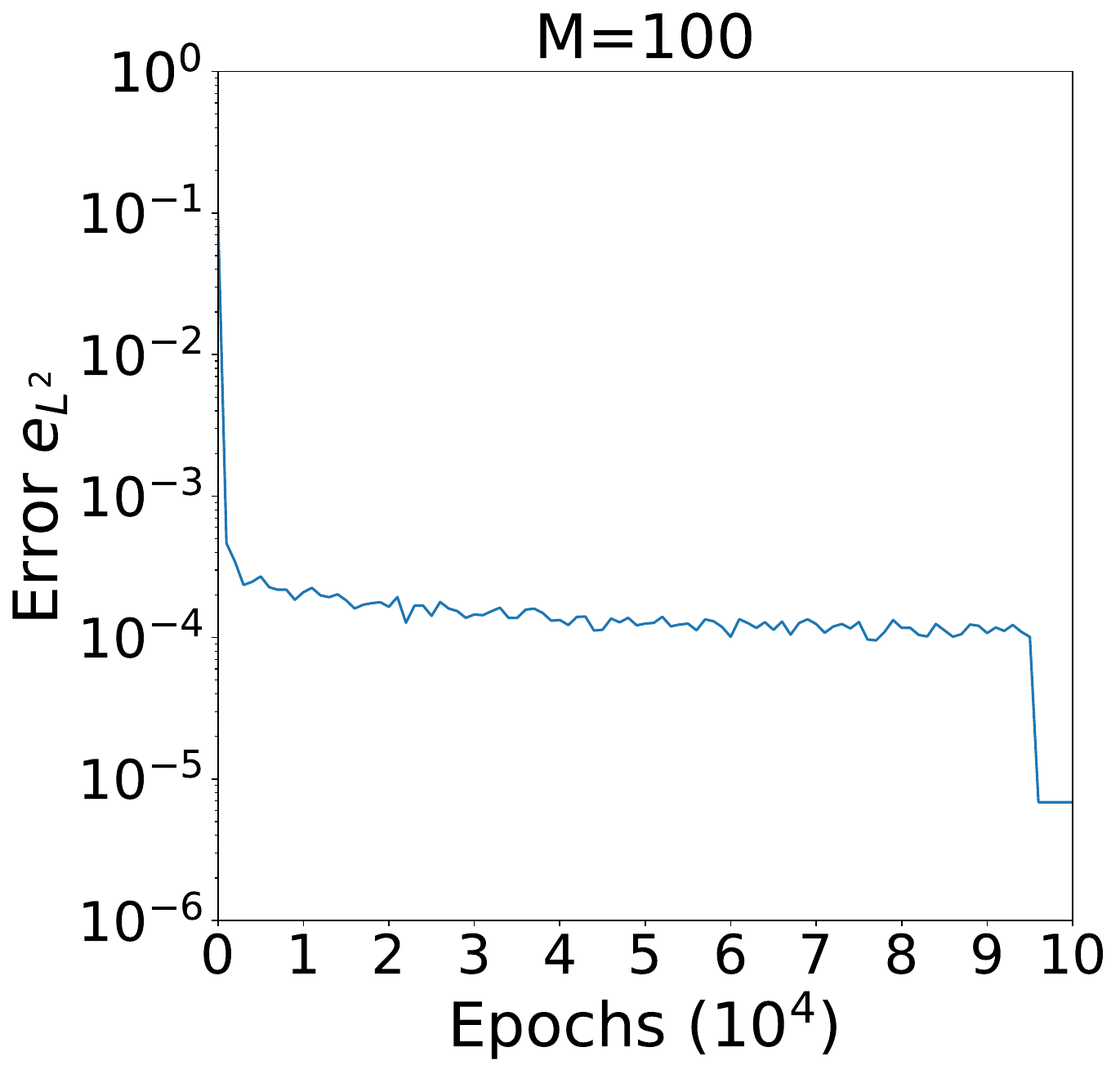}
\includegraphics[width=3.5cm,height=3.5cm]{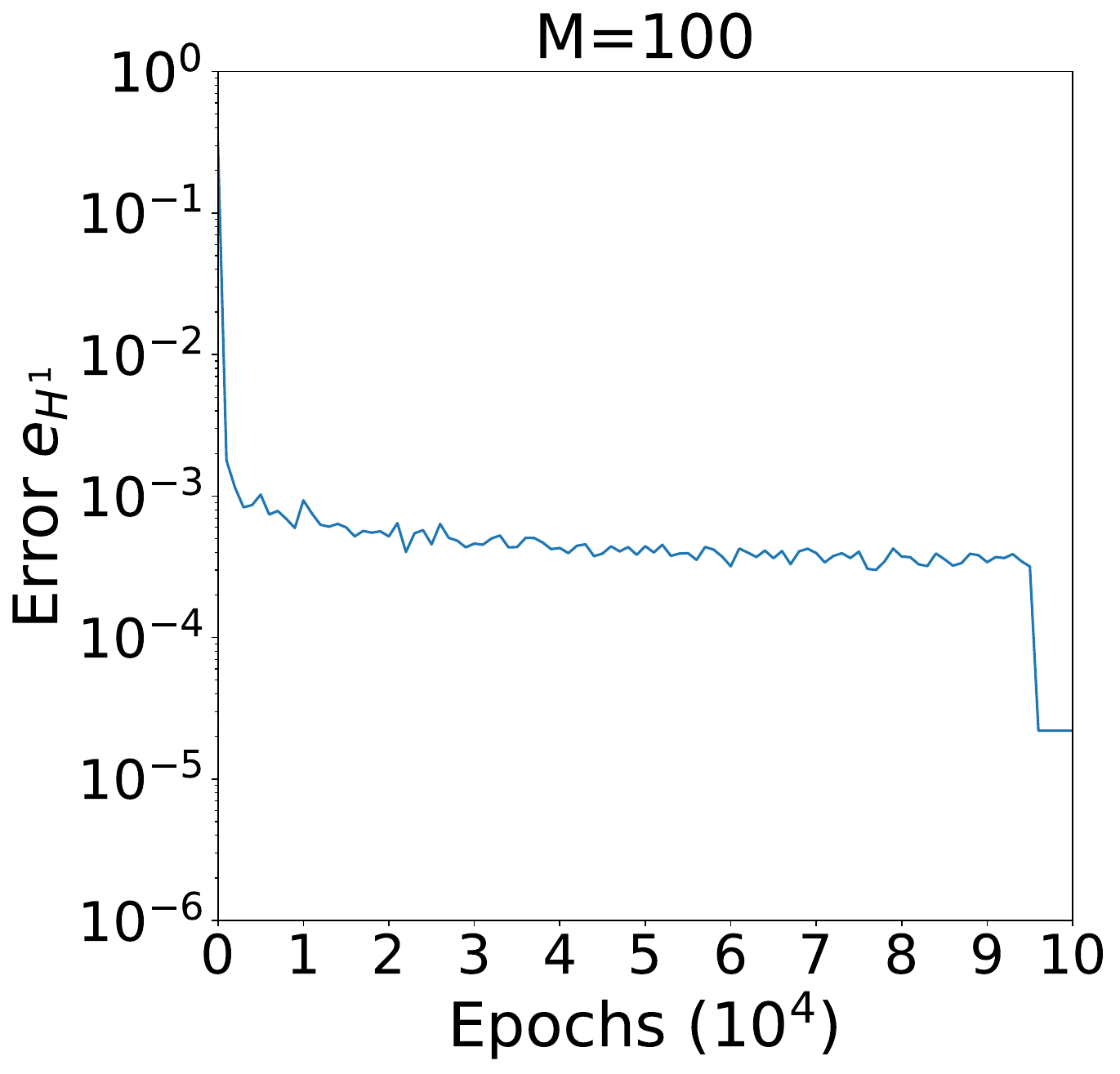}
\caption{The relative errors during training process in Example 3.}\label{fig_errors_ex3}
\end{figure}

The final errors are listed in the 
Table \ref{table_errors_ex3} lists the final errors, which shows that 
the TNN based machine learning method can solve the high-dimensional stochastic
partial differential equations with high accuracy.   
The relative errors $e_{L^2}$ and $e_{H^1}$ versus the number of 
epochs are shown in Figure \ref{fig_errors_ex3} .
\begin{table}[htb!]
\caption{Errors of Example 3.}\label{table_errors_ex3}
\begin{center}
\begin{tabular}{ccccc}
\hline
$M$ &   $e_{L^2}$ & $e_{H^1}$\\
\hline
10  &   3.609e-08 & 1.141e-07 \\
20  &   5.958e-07 & 1.872e-06 \\
50  &   5.983e-06 & 1.898e-05 \\
100 &   6.855e-06 & 2.198e-05 \\
\hline
\end{tabular}
\end{center}
\end{table}

\section{Conclusion}
In this paper, we design a type of TNN based machine learning method for
solving the elliptic partial differential equations with random coefficients
in a bounded physical domain.
Different from the general FNN-based machine learning method,
TNN has the tensor product structure and then the corresponding 
high-dimensional integration can be computed with high accuracy.
Benefit from the high accuracy of the high-dimensional integration,
the TNN based machine learning method can solve the high-dimensional
differential equations with high accuracy.
We believe that the ability of TNN based machine learning method will bring
more applications in solving linear and nonlinear high-dimensional PDEs. 
These will be our future work.

\end{document}